\documentclass{article}
\usepackage[utf8]{inputenc}
\usepackage[T1]{fontenc}
\usepackage{amsmath,amsfonts,amssymb,amsthm}
\usepackage[a4paper,left=2.cm,right=2.cm,top=2.cm,bottom=2.cm]{geometry}
\usepackage{algorithm,algpseudocode}
\usepackage{color}
\usepackage{graphicx}
\usepackage{svg}
\usepackage{array}
\usepackage{subcaption}
\usepackage{pgf}
\newcolumntype{C}[1]{>{\centering\arraybackslash}m{#1}}

\theoremstyle{remark}
\newtheorem{remark}{Remark}

\newcommand{\R}{\ensuremath{\mathbb{R}}}
\newcommand{\bx}{\ensuremath{\mathbf{x}}}
\newcommand{\by}{\ensuremath{\mathbf{y}}}
\newcommand{\br}{\ensuremath{\mathbf{r}}}
\newcommand{\bz}{\ensuremath{\mathbf{z}}}
\newcommand{\bb}{\ensuremath{\mathbf{b}}}
\newcommand{\bq}{\ensuremath{\mathbf{q}}}

\newcommand{\bu}{\ensuremath{\mathbf{u}}}

\newcommand{\bw}{\ensuremath{\mathbf{w}}}

\newcommand{\bv}{\ensuremath{\mathbf{v}}}
\newcommand{\bV}{\ensuremath{\mathbf{V}}}
\newcommand{\bT}{\ensuremath{\mathbf{T}}}
\newcommand{\bL}{\ensuremath{\mathbf{L}}}
\newcommand{\bX}{\ensuremath{\mathbf{X}}}
\newcommand{\bA}{\ensuremath{\mathbf{A}}}
\newcommand{\bM}{\ensuremath{\mathbf{M}}}

\newcommand{\bZ}{\ensuremath{\mathbf{Z}}}
\newcommand{\bC}{\ensuremath{\mathbf{C}}}
\newcommand{\bI}{\ensuremath{\mathbf{I}}}
\newcommand{\bJ}{\ensuremath{\mathbf{J}}}
\newcommand{\bW}{\ensuremath{\mathbf{W}}}
\newcommand{\bP}{\ensuremath{\mathbf{P}}}
\newcommand{\bR}{\ensuremath{\mathbf{R}}}
\newcommand{\bQ}{\ensuremath{\mathbf{Q}}}
\newcommand{\bS}{\ensuremath{\mathbf{S}}}
\newcommand{\bdelta}{\ensuremath{\boldsymbol{\delta}}}

\newcommand{\bTheta}{\ensuremath{\boldsymbol{\Theta}}}
\newcommand{\bDel}{\ensuremath{\boldsymbol{\Delta}}}
\newcommand{\bXi}{\ensuremath{\boldsymbol{\Xi}}}

\newcommand{\bone}{\ensuremath{\boldsymbol{1}}}

\newcommand{\sK}{\ensuremath{\mathcal{K}}}

\newcommand{\tbx}{\ensuremath{\mathbf{\tilde{x}}}}
\newcommand{\hbz}{\ensuremath{\mathbf{\hat{z}}}}
\newcommand{\hbZ}{\ensuremath{\mathbf{\hat{Z}}}}

\author{Ahmed Chabib, Jean-François Witz, Vincent Magnier, Pierre Gosselet\\
Univ. Lille, CNRS, Centrale Lille, UMR 9013 - LaMcube - \\ Laboratoire de Mécanique, Multiphysique, Multiéchelle, \\ F-59000 Lille, France
}

\title{Conjugate gradient for ill-posed problems: regularization by preconditioning, preconditioning by regularization}

\date{\today}
\begin{document}
	\maketitle

\begin{abstract}This paper investigates using the conjugate gradient iterative solver for ill-posed problems. We show that preconditioner and Tikhonov-regularization work in conjunction. In particular when they employ the same symmetric positive semi-definite operator, a powerful Ritz analysis allows one to estimate at negligible computational cost the solution for any Tikhonov's weight. This enhanced linear solver is applied to the boundary data completion problem and as the inner solver for the optical flow estimator.
\end{abstract}

\section{Introduction}

Ill-posed systems of equations are ubiquitous in mechanics. They are particularly present in identification problems, such as the boundary completion in elasticity \cite{2011_Kadri_Crack,2018_Ferrier_Ritz}. They also appear in methods involving some compact operator, like the Herglotz' transform to build solutions to the Helmholtz problems \cite{2015_Kovalesky_complex_ray}. Beside issues of existence and uniqueness, ill-posed problems are characterized by the lack of stability between the cause and the effect, in other words small perturbations in the input potentially lead to large modifications of the output. 

In this paper, we focus on discrete $n\times n$ linear(ized) symmetric positive semi-definite systems of the form $\bA\bx=\bb$, allowing to analyze all properties in terms of the spectrum of $\bA$ which can be diagonalized as $\bA=\sum_{i=1}^n \sigma_i \bu_i\bu_i^T$ with $(\bu_i)$ an orthonormal basis and $(\sigma_i)$ non-negative eigenvalues in decreasing order. 

Existence and uniqueness are linked to the null-space of $\bA$ (strictly zeros eigenvalues) whereas stability is associated with the accumulation of eigenvalues near zero. Indeed, a small contribution of $\bb$ in an eigendirection associated with a small eigenvalue of $\bA$ has a significant impact on $\bx=\sum \frac{\bu_i^T\bb}{\sigma_i}\bu_i$. Ill-posed problems thus result in poorly conditioned operators.

Solving such systems amounts to finding a satisfactory treatment to these small eigenvalues: truncation, shift, filtering.
\begin{itemize}
	\item  Truncation involves disregarding the problematic directions, and only keeping the part of the matrix associated with eigenvalues larger than a given criterion~$\varepsilon$: 
	\begin{equation}
		\bA^{-1}\simeq \sum_{i=1}^m \frac{\bu_i\bu_i^T}{\sigma_i} ,\qquad \sigma_i>\varepsilon>0\text{ for } i\leqslant m<n.
	\end{equation}
	This method easily extends to general matrices using the singular value decomposition \cite{1987_Hansen_tsvd}.
	\item Shift is generally achieved thanks to Tikhonov regularization \cite{1977_Tikhonov_solutions}, that can be written as:
	\begin{equation}
		\bA \simeq  \bA_\lambda = \bA+\lambda\bM.
	\end{equation}
	The simplest choice is $\bM=\bI$ and in that case $\lambda>0$ becomes the lower bound of the spectrum of $\bA_\lambda$. It is of course preferable to use a matrix $\bM$ with more physical sense, acting more locally on the small eigenvalues. Indeed, a ``flat'' regularization like $\lambda\bI$ may temper the contribution of the higher part of the spectrum. 	
	For instance one could use $\lambda(\sum_{i>m}\bu_i\bu_i^T)$, and the shift would behave like the truncation when $\lambda\to\infty$.
	
	\item Filtering tries either to ``clean'' the right-hand side from spurious excitation or to improve the solution after it was computed by enforcing some physical properties. For instance, smoothing can be used to recover regularity in a noisy solution. These approaches can be implemented by projectors: the former by right-projection $\bx = \bA^{-1}\bP_r \bb$, and the latter by left projection $\bx = \bP_l\bA^{-1}\bb$. Note that the approaches are equivalent when the preconditioner is orthogonal on an eigensubspace, $\bP = \bP^T = \sum_{i\leqslant m} \bu_i\bu_i^T$ because in that case $\bP\bA=\bA\bP$ and there is a direct correspondence between regular excitation and regular solution. Another example of symmetric filtering is the use of coarse discretization in Galerkin methods. For instance, in Digital Image Correlation \cite{2006_Besnard_FE}, the coarser the discretization the better posed the problem -- but of course the less precise.
\end{itemize}

All these techniques are often controlled by a parameter ($\varepsilon$ for the truncation, $\lambda$ for the regularization,\ldots) which needs to be tuned in order to find a balance between the information inside the original system and the information brought (or removed) by the treatment.
When the accuracy of the data is known, Morozov's principle \cite{1968_Morozov_principle} provides an objective criterion for choosing the parameter: the correction introduced by the added information should not exceed the noise in the measurement.

When no such data is available, a compromise must be found. Picard's principle \cite{1990_Hansen_Picard} compares the eigenvalues $(\sigma_i)$ (sorted in decreasing order) and the decomposition of the right-hand side on the eigendirections $(\bu_i^T\bb)$. While eigenvalues decrease less rapidly than their contributions to the right-hand side, the solution remains controlled. The L-curve \cite{1992_Hansen_Lcurve} is a visual aid to find a balance. The solutions for various levels of regularization are positioned  in a frame (``norm of the residual'', ``norm of the solution''). In general large regularization leads to low norm of the solution but high error, whereas small regularization leads to lower level of error but large solutions (highly perturbed). Ideally, some corner exists which realizes a trade-off between residual and oscillating solution.\medskip

Solving an ill-posed system with an iterative solver may seem counterintuitive, as poor preconditioning is often considered a red flag for the use of such solvers. In fact, it can be viewed as an opportunity to implement approximately the strategies mentioned above: because of their kinship with power iteration, iterative solvers have the tendency to favor the higher part of the spectrum in the beginning, and the limited number of iterations (often controlled by the convergence threshold) can be a way to stop iterations before searching too deeply in the ill-conditioned portion of the spectrum. This paper focuses on the many opportunities opened by the introduction of a preconditioner in the iteration. 

Preconditioner $\bM$ is often introduced in numerical methods courses as a cheap-to-compute approximation of the inverse of the problem's operator ($\bM^{-1}\simeq\bA^{-1}$) in the sense that the spectrum of $\bM^{-1}\bA$ should be as concentrated as possible around a non-zero value (which can be scaled to $1$). This can be roughly estimated by the condition number of $\bM^{-1}\bA$, but more sophisticated studies are available \cite{AXELSSON:1986:RCPCG}. However this does not apply for ill-posed problems where the inverse of $\bA$ can not be properly defined. Anyhow, the preconditioner may play a major role as it can help tone the bad part of the spectrum down. It is well-known \cite{SAAD:1992:NMLEP} that preconditioning with SPD matrix $\bM$ is equivalent to solving the system $\bL^{-1}\bA\bL^{-T}\by = \bL^{-1}\bb$ where $\bL\bL^T=\bM$ is the Cholesky factorization, and $\bx = \bL^{-T}\by$. This system is governed by the eigenvalues of $\bL^{-1}\bA\bL^{-T}$ or more practically by the generalized eigenvalues $(\mu_i)$ of $(\bA,\bM)$: $\bA\bv_i=\mu_i\bM\bv_i$, with $(\bv_i)$ a $\bM$-orthonormal basis. Using Rayleigh quotient $R(\bv)=(\bv^T\bA\bv)/(\bv^T\bM\bv)$, one clearly sees that if the denominator is large enough in some directions $(\bv)$, the associated eigenvalues will be reduced and this part of the search space will be explored later in the iterations. We see that the choice of the preconditioner obeys similar criteria as the choice of the regularization, hence the use of the same notation $\bM$. 

In the case of well-posed problems, a precise solution is achievable whatever the preconditioner which mainly impacts the number of iterations. In the case of ill-posed problem, preconditioning can be a tool to favor ``better'' (more regular) directions. Since trying to achieve a too strict convergence may not be realistic, the sequence of directions suggested by the preconditioner plays an important role in the definition of retained solution for a relatively weak convergence threshold.

The use of preconditioned conjugate gradient seems not to be that frequent when solving ill-posed problems, other Krylov solvers  based of least-square (MinRes, CGNE, CGLS)~\cite{CGminres12,WANG2024147,Neubauer2022,Buccini2017} or Landweber iteration~\cite{XIONG2017108} are often preferred due to certain monotony properties. We believe that our study offer a set of important tools to reconsider the potential of conjugate gradient for ill-posed problems. 

\subsection*{Main contributions} In this paper, we attempt to combine these ideas within a sophisticated preconditioned conjugate gradient solver with several contributions. Note that we do not pretend to propose a new regularization method, but we think that we bring new insight on the role of preconditioning and a particularly well-tuned solver where all regularization tools can be deployed at marginal cost.

First we show that:
\begin{itemize}
	\item The preconditioner can be leveraged to naturally regularize the solution.
	\item The preconditioner provides appropriate norms to analyze the properties of the solution and plot exploitable L-curves.
	\item In the spirit of~\cite{2018_Ferrier_Ritz} we show that the Ritz analysis may allow an interesting a posteriori filtering of the solution. 
\end{itemize}

A key novelty lies in the cases where a Tikhonov regularization with an easily solvable structure is used, where we show that it is possible:
\begin{itemize}
	\item To fully analyze the effect of the regularization on the original system in particular through Picard plots.
	\item To post-process at negligible cost an approximation of the solution for any other regularization weight $\lambda$. 
\end{itemize}

When embedded in a nonlinear, process we show that:
\begin{itemize}
\item Acceleration is available.
\item It is possible to postprocess approximations of the solutions for a predetermined family of weights $(\lambda_i)$ at negligible cost.
\end{itemize}  

For our assessments, we consider the Steklov-Poincaré approach to complete boundary data which has the advantage of being a linear problem, and the recovery of the optical flow between two images, which is nonlinear. 

\subsection*{Organization of the paper}
The paper is organized as follows.  In Section~\ref{sec:PCG}, we recall the augmented preconditioned conjugate gradient algorithm and the computation of Ritz eigenelements, we discuss and illustrate the effect of preconditioning and regularizing separately and the norms to analyze the evolution of the iteration. Particularly, in Section~\ref{sec:PrecReg}, we consider the case of regularized systems preconditioned by the regularization matrix where we can analyze in detail the effect of the regularization. These properties are assessed on the data completion problem in Section~\ref{sec:datacompletion}. Section~\ref{sec:APCG} presents the
adjustments required to handle nonlinear problems, in particular augmentation.  Section~\ref{sec:OF} presents the assessments for estimating the optical flow. Section~\ref{sec:conc} concludes the paper.

\subsection*{Notations}
We use normal font for scalars, boldface lowercase for vectors and boldface uppercase for matrices. A collection of vectors $(\bx_j)$ can be put in the matrix form $\bX_m=(\bx_0,\ldots,\bx_{m-1})$, the index $m$ thus corresponds to the number of columns of the matrix. This work is presented in $\mathbb{R}^n$ even though the methods also apply for complex Hermitian matrices and vectors.

\section{Preconditioned Conjugate Gradient and Ritz elements}\label{sec:PCG}

Let $\bA$ be a symmetric definite positive matrix and $\bb$ be a vector. We search the solution to the system $\bA\bx=\bb$. We use a conjugate gradient, preconditioned by the symmetric positive semi-definite matrix $\bM$.

At iteration $i$, we note $\bx_i$ the approximation and $\br_i=\bb-\bA\bx_i$ the residual. We introduce the Krylov subspace
$\sK_i(\bM^{-1}\bA,\bM^{-1}\br_0)$ \cite{SAAD:1997:DAK}:
\begin{equation}
	\sK_i(\bM^{-1}\bA,\bM^{-1}\br_0) = \\
	\operatorname{span}\left(\bM^{-1}\br_0,\ldots,(\bM^{-1}\bA)^{(i-1)}\bM^{-1}\br_0\right)
\end{equation}

Given an arbitrary initialization $\bx_{0}$ and associated residual $\br_{0}=\bb-\bA\bx_{0}$, the $i_{th}$ iteration can be defined as:
\begin{equation}\label{eq:krysol_principle}
	\left\{\begin{array}{ll}
		\textrm{find } &\bx_i \in \bx_{0}+\sK_i(\bM^{-1}\bA,\bM^{-1}\br_{0})\\
		\textrm{such that }&\br_i \perp \sK_i(\bM^{-1}\bA,\bM^{-1}\br_{0})
	\end{array}\right.
\end{equation}
This iteration is achieved by Algorithm~\ref{alg:PCG}. 

\begin{algorithm}[ht]
	\caption{Preconditioned Conjugate Gradient}
	\label{alg:PCG}
	\begin{algorithmic}
		\State $\br_0=  \bb - \bA \bx_{0}=\bP^T\br_{0}$
		\State {${\bz_{0}} = \bP\bM^{-1}   \br_{0}$, $\bw_0=\bz_0$}
		\State $\gamma_0 = (\bz_{0}^T \br_{0})$
		\For{$i = 0,\,1,\, \dots,m$ (convergence)}
		\State  $\bq_i = \bA{\bw_{i}}$
		\State  $\delta_i = (\bw_i^T \bq_{i})$,   $\alpha_i = \delta_i^{-1}\gamma_i$ 
		\State   $\bx_{i+1} = \bx_{i}+\bw_i \alpha_i $
		\State   $\br_{i+1} = \br_i- \bq_i \alpha_i$
		\State   {${\bz_{i+1}} =\bP\bM^{-1} \br_{i+1} $}
		\State $\gamma_{i+1} = (\bz_{i+1}^T \br_{i+1})$
		\State  $\beta_{i}= \gamma_i^{-1}\gamma_{i+1}$
		\State $\bw_{i+1} = \bz_{i+1} + \bw_i\beta_{i} $
		\EndFor
	\end{algorithmic}
\end{algorithm}
The algorithm builds two special basis of $\sK_i(\bM^{-1}\bA,\bM^{-1}\br_0)$, $\bZ_i$ is $\bM$-orthogonal whereas $\bW_i$ is $\bA$-orthogonal: 
\begin{equation}
	\begin{aligned}
		\bZ_i^T\bM\bZ_i = \bZ_i^T\bR_i=\operatorname{diag}(\gamma_j)_{0\leqslant j<i} \\
		\bW_i^T\bA\bW_i = \bW_i^T\bQ_i=\operatorname{diag}(\delta_j)_{0\leqslant j<i} \\
	\end{aligned}
\end{equation}

It is convenient to introduce the $\bM$-normalized version of the $\bZ_i$ basis:
\begin{equation}
	\hbz_i = \frac{(-1)^i\bz_i}{\sqrt{\gamma_i}}\qquad 
	\text{so that } \hbZ_i^T\bM\hbZ_i = \bI 
\end{equation}	
$\hbZ_i$ is in fact the basis that would have been obtained by the Arnoldi procedure \cite{2003_Saad_iterative}, and we have:	
\begin{equation}
	\begin{aligned}
		\hbZ_i^T\bA\hbZ_i &= \bT_i = \operatorname{Tridiag}(\eta_{j-1},\mu_{j},\eta_{j}) \\
		&\text{with } \mu_0 = \frac{1}{\alpha_0},
		\quad \mu_j=\frac{1}{\alpha_j}+\frac{\beta_{j-1}}{\alpha_{j-1}},\quad \eta_{j}= \frac{\sqrt{\beta_j}}{\alpha_j}
	\end{aligned}
\end{equation}
We can diagonalize $\bT_i=\bXi_i\bTheta_i\bXi_i^T$ where $\bTheta_i$ is the diagonal matrix of eigenvalues sorted in decreasing order and $\bXi_i$ the orthonormal matrix of eigenvectors. 

The Ritz vectors are $\bV_i=\hbZ_i\bXi_i$, while $\bTheta_i$ are the Ritz values of the system. They satisfy:
\begin{equation}
	\bV_i^T \bM \bV_i=\bI \qquad\text{and}\qquad \bV_i^T \bA \bV_i=\bTheta_i.
\end{equation}
In order to mark the dependency of the Ritz vectors and values on the iteration $i$, they are denoted with an exponent $(i)$: $\bTheta_i=\operatorname{diag}(\theta_j^{(i)})_{1\leqslant j\leqslant i}$ and  $\bV_i=\begin{pmatrix}
	\bv_1^{(i)},\ldots,\bv^{(i)}_i
\end{pmatrix}$. As the number of iterations $i$ increases, the $(\theta^{(i)}_{j})_{1\leqslant j\leqslant i}$ and $(\bv^{(i)}_{j})_{1\leqslant j\leqslant i}$ tend to approximate the generalized eigenvalues and eigenvectors of the couple $(\bA,\bM)$ \cite{JIA:2004:CRV}.

\subsection{Appropriate solution and error norms and stopping criteria}

Conjugate gradient gives valuable pieces of information at no cost in the course of the iterations, but in specific norms. Indeed, the preconditioner can be viewed as providing a physic-based alternative to the simple Euclidean norm and inner product.

First, we have error estimators~\cite{hestenes1952,2001_Axelsson_norm_PCG}:
\begin{equation}
	\begin{aligned}
		\|\br_i\|_{\bM^{-1}}^2 &=\gamma_i \\ 
		\|\bx_{i+1}-\bx\|_\bA^2 &= \|\bx_i-\bx\|_\bA^2 - \gamma_i^2\delta_i
	\end{aligned}
\end{equation}
of course the difficulty for the second identity is that $\|\bx_0-\bx\|_\bA^2$ is unknown. Note that in~\cite{Meuran2021}, a strategy is proposed to overcome this problem and to fully estimate  $\|\bx_{i+1}-\bx\|_\bA^2$.

We also have measurement of the norm of the correction brought by iterations \cite{2018_Ferrier_Ritz}:
\begin{equation}
	\begin{aligned}
		\|\bx_{i+1}-\bx_0\|_{\bM}^2 &= \|\bx_{i}-\bx_0\|_{\bM}^2 + \alpha_i^2\|\bw_i\|^2_{\bM}+2\alpha_i(\bw_i^T\bM(\bx_i-\bx_0))\\
		\text{with } &\left\{\begin{aligned}
			&\|\bw_{i+1}\|^2_{\bM} = \gamma_i+\beta_{i}^2\|\bw_i\|^2_{\bM}, \qquad \|\bw_0\|^2_{\bM}=\gamma_0,\\
			&(\bw_{i+1}^T\bM(\bx_{i+1}-\bx_0))=-\beta_{i}\left((\bw_i^T\bM(\bx_i-\bx_0))+\alpha_i\|\bw_i\|_{\bM}^2\right).
		\end{aligned}\right.
	\end{aligned}
\end{equation}
Finally, we have an estimator on the preconditioned operator:
\begin{equation}
	\begin{aligned}
		\|\bT_{0}\|_F^2 &= \mu_0^2,\\
		\|\bT_{i+1}\|_F^2 &=  \|\bT_i\|_F^2 + \mu_i^2 + \eta_i^2 +\eta_{i-1}^2 \to \|\bM^{-\frac{1}{2}}\bA\bM^{-\frac{1}{2}}\|^2_{F}, 
	\end{aligned}
\end{equation}
where index $F$ stands for the Frobenius norm, $\|\bM^{-\frac{1}{2}}\bA\bM^{-\frac{1}{2}}\|^2_{F}$ is the sum of the squares of the generalized eigenvalues of $(\bA,\bM)$.\medskip

We can then devise costless (without extra computation) stopping criteria: 
\begin{subequations}\label{eq:stopping}
	\begin{equation}
		\|\br_i\|_{\bM^{-1}} < \varepsilon \|\br_0\|_{\bM^{-1}}, \label{eq:stopping1}
			\end{equation}
			\begin{equation}
		\|\br_i\|_{\bM^{-1}} < \varepsilon \|\bT_i\|_{F} \|\bx_i-\bx_0\|_{\bM},\label{eq:stopping2}
		\end{equation}
		\begin{equation}
		\gamma_i^2\delta_i^{-1} < \varepsilon^2. \label{eq:stopping3}
	\end{equation}
\end{subequations}
The first one is very classical, but it is risky in the sense that it may be too strict if the initialization was well-chosen ($\|\br_0\|_{\bM^{-1}}$ is already small). The second one is inspired from the Scipy implementation of MinRes with a more adapted choice of norms, we were unable to trace the original source of this idea. This criterion is useful because it balances the error reduction against the growth of the solution norm, a central dilemma when solving ill-posed problems. The third criterion corresponds to the iteration bringing negligible correction, what is often called stagnation, in general it must be checked for a sequence of iterations before actually stopping. 
It is often interesting to combine the criteria, add stagnation detection, and to also use safeguards in absolute value in case of too good initialization.

\subsection{Natural frame for the L-curve during iterations}
In the case of a poorly conditioned system, the reduction of the error can be obtained at the price of an explosion of the norm of the solution. This is well explained by Picard analysis: the phenomenon occurs when the eigenvalues of the operator decrease faster than the contribution of the right-hand side in the associated direction. It can also be visualized on a L-curve, in the positive quarter of a frame of the form $(\|\br_i\|, \|\bx_{i}\|)$: the curve starts in the bottom right corner (large error, small norm) with a fast decay of the error, and finishes in the top left corner (reduced error, large norm). As shown earlier, conjugate gradient provides natural norms to evaluate the error and the norm of the solution: $\|\bx_i-\bx\|_\bA$ and $\|\bx_i-\bx_0\|_\bM$. With this choice of norms, the curve is always oriented toward the upper-left corner: at each iteration, the norm of the error decreases and the norm of the solution increases.

\subsection{A posteriori filtrering}
Ritz elements offer a convenient way to filter the solution. Assuming $m$ iterations were conducted, we can process the basis $\bV_m$ and the values $\bTheta_m$. We can decompose the right-hand side on the Ritz basis $r^{(m)}_j=\bv_j^{{(m)}^T}\br_0$, and define:
\begin{equation}
	\text{for }i\leqslant m,\qquad 	\tbx^{(m)}_i = \bx_0 + \sum_{j=1}^i  \frac{r^{(m)}_j}{\theta^{(m)}_j}\bv^{(m)}_j.
\end{equation}
We have:
\begin{equation}
	\begin{aligned}
		\|\tbx^{(m)}_i-\bx\|_\bA^2 &= \|\tbx^{(m)}_i-\bx_0\|_\bA^2 - \sum_{j=1}^i \frac{(r^{(m)}_j)^2}{\theta^{(m)}_j},\\
		\|\tbx^{(m)}_i-\bx_0\|_\bM^2 &=  \sum_{j=1}^i \frac{(r^{(m)}_j)^2}{(\theta^{(m)}_j)^2},
	\end{aligned}
\end{equation}
and of course:
\begin{equation}
	\|\tbx^{(m)}_i-\tbx^{(m)}_{i-1}\|_\bA^2 =   \frac{(r^{(m)}_i)^2}{\theta^{(m)}_i} \qquad \text{and}\qquad \|\tbx^{(m)}_i-\tbx^{(m)}_{i-1}\|_\bM^2 =   \frac{(r^{(m)}_i)^2}{(\theta_i^{(m)})^2}.
\end{equation}
Since the $(\theta^{(m)}_i)$ are sorted in decreasing order, we see that the error of $i\mapsto(\tbx^{(m)}_i)$ tends to decrease slower than its norm tends to increase. The L-curve for $(\tbx^{(m)}_i)_i$ is then convex and does not zig-zag. 

The slope of the L-curve between the point $i-1$ and $i$ is $-(\theta^{(m)}_i)^{-1}$. A possibility is to define the corner as the point which maximizes the variation of slope:  $i=\arg\max_j((\theta^{(m)}_{j+1})^{-1}-(\theta^{(m)}_j)^{-1})$. 

Ritz' elements also make it possible to use Picard's theory and stop the construction of $\tbx^{(m)}_i$ when the contribution $j\mapsto r^{(m)}_j$ starts to decrease less fast than $j\mapsto\theta^{(m)}_j$. This criterion has the advantage to take into account the properties of the right-hand side.

\section{Preconditioning by regularization}\label{sec:PrecReg}


The core idea of this paper is that the preconditioner and the regularization should be the same operator. Thus, we use the matrix $\bM$ for the Tikhonov regularization. As mentioned in the introduction, this idea makes sense as the same physical motivation underlies the choice of the regularization and that of the preconditioner. Moreover, many opportunities are opened by this choice. However, the hypothesis behind this idea is that there exists a cheap technique to apply the preconditioner (i.e. $\bM^{-1}$).

We are interested in Tikhonov-regularized systems of the form:
\begin{equation}\label{eq:regusys}
	\underset{\bA_\lambda}{\underbrace{(\bA + \lambda \bM )}} \bx_\lambda = \underset{\bb_\lambda }{\underbrace { \bb_{\bA} + \lambda \bb_{\bM}}}.
\end{equation}
Note that we chose a $\lambda$-affine form for the right-hand side because it meets our practical needs, but the method applies to any separate form ($\bb(\lambda) = \sum_a f_a(\lambda)\bb_a$).

If we assume that the system~\eqref{eq:regusys} was solved for a given $\lambda$ in $m$ iterations using the proposed $\bM$-preconditioned conjugate gradient algorithm, then we can process the Ritz basis $\bV_m$. The strong point is that $\bV_m$ separates the effects of the operator $\bA$ and that of the regularization $\bM$ independently of $\lambda$:
\begin{equation}\label{eq:specritz}
	\begin{aligned}
		\bV_m^T\bM\bV_m &= \bI_m, \\
		\bV_m^T\bA_\lambda\bV_m &= \bTheta_{\lambda,m}=\bTheta_m+\lambda\bI_m.
	\end{aligned}
\end{equation}
\begin{remark}
	$\lambda$ can be viewed as a shift in the generalized eigenvalues of $(\bA,\bM)$. Since $\lambda$ alters the initial residual and only a limited number $m$ of iterations is made, the content of $\bV_m$ depends on the choice of $\lambda$, but the orthogonality properties remain.
\end{remark}

Note that the initial residual takes the separate form:
\begin{equation}
	\br_{\lambda,0}=\bb_\lambda-\bA_\lambda\bx_0 = \underbrace{(\bb_{\bA} - \bA\bx_0)}_{\br_{\bA,0}}+ \lambda\underbrace{(\bb_{\bM} - \bM\bx_0)}_{\br_{\bM,0}}, 
\end{equation}
and we can define the spectral contributions $r^{(m)}_{A,j}=\bv_j^{{(m)}^T}\br_{\bA,0}$ and $r^{(m)}_{M,j}=\bv_j^{{(m)}^T}\br_{\bM,0}$,

After $m$ iterations, we can define the Ritz' approximations for $i\leqslant m$:
\begin{equation}\label{eq:ritzapprox}
	\tbx^{(m)}_{\lambda,i} = \bx_{0} + \sum_{j=1}^i \bv^{(m)}_j \frac{(\bv_j^{(m)^T}\br_{\lambda,0})}{\theta^{(m)}_j+\lambda} = \bx_{0} + \sum_{j=1}^i  \frac{r^{(m)}_{A,j}+\lambda r^{(m)}_{M,j}}{\theta^{(m)}_j+\lambda}\bv^{(m)}_j
\end{equation}
These approximations can be computed at marginal cost, and there dependence in $\lambda$ is explicit: the L-curve of $\lambda\mapsto \tbx^{(m)}_{\lambda,i}$ is a rational fraction. It even permits to give sense to the limit solution when $\lambda \to 0$ even when $\bA$ was not invertible. It also gives an analytical formula for the search of the optimal choice of $(\lambda,i)$ realizing a good compromise between error and norm of the solution.
 we have the properties:
\begin{equation}\label{eq:l-curves}
	\begin{aligned}
		\|\tilde{\bx}^{(m)}_{\lambda,i} - \bx_0\|^2_\bM  &= \sum_{j=1}^i \left(\frac{r^{(m)}_{A,j}+\lambda r^{(m)}_{M,j}}{\theta^{(m)}_j+\lambda}\right)^2\\ 
		\|\tilde{\bx}^{(m)}_{\lambda,i} - \bx_\lambda\|^2_{\bA_\lambda}  &= \|\tilde{\bx}^{(m)}_{\lambda,0} - \bx_\lambda\|^2_{\bA_\lambda}-\sum_{j=1}^i\frac{\left(r^{(m)}_{A,j}+\lambda r^{(m)}_{M,j}\right)^2}{\theta^{(m)}_j+\lambda}
	\end{aligned}
\end{equation}

We write $\bV^{(m)}_i=\begin{pmatrix}
	\bv^{(m)}_1,\ldots,\bv^{(m)}_i
\end{pmatrix}$ and $\bTheta^{(m)}_i=\operatorname{diag}(\theta_1^{(m)},\ldots,\theta_i^{(m)})$, so that $\tilde{\bx}^{(m)}_{\lambda,i}= \bx_0+\bV^{(m)}_i (\bTheta^{(m)}_i+\lambda
\bI_i)^{-1}\bV^{(m)^T}_i\br_{0,\lambda}$.
We have an even more interesting measure of the error in non-regularized norm with respect to the non-regularized solution:
\begin{equation}\label{eq:l-curves2}
	\begin{aligned}
		\|\tilde{\bx}^{(m)}_{\lambda,i} - \bx\|^2_{\bA} - \|\bx_0- \bx\|^2_{\bA} &= \|\bV^{(m)}_i (\bTheta^{(m)}_i+\lambda
		\bI_i)^{-1}\bV^{(m)^T}_i\br_{\lambda,0} + \bx_0- \bx\|^2_{\bA} - \|\bx_0- \bx\|^2_{\bA} \\
		&=  \|\bV^{(m)}_i (\bTheta^{(m)}_i+\lambda
		\bI_i)^{-1}\bV^{(m)^T}_i\br_{\lambda,0}\|^2_{\bA} \\
		&\qquad\qquad- 2\br_{\bA,0}^T \bV^{(m)}_i (\bTheta^{(m)}_i+\lambda
		\bI_i)^{-1}\bV^{(m)^T}_i\br_{\lambda,0} \\
		&= \sum_{j=1}^i  \frac{\left(r^{(m)}_{A,j}+\lambda r^{(m)}_{M,j}\right)}{(\theta^{(m)}_j+\lambda)}\left( \frac{\theta^{(m)}_j\left(r^{(m)}_{A,j}+\lambda r^{(m)}_{M,j}\right)} {\left(\theta^{(m)}_j+\lambda\right)} - 2 r^{(m)}_{A,j}\right)
	\end{aligned}
\end{equation}

Another feature is the possibility to analyze the system in terms of spectral content ($\theta^{(m)}_i$) vs regularization $\lambda$, and with regard to the decomposition of the right-hand side $\left(r^{(m)}_{A,j}+\lambda r^{(m)}_{M,j}\right)$. This can be particularly well visualized in a Picard plot.

\section{Assessment in the linear case: data completion problem}\label{sec:datacompletion}
\subsection{Laplace PDE with missing and redundant boundaries}
We use the classical illustration of the ill-posedness for the inverse Laplace problem. We consider the rectangular domain $\Omega:=[0,T]\times[0,H]$, where the following Cauchy problem holds:
\begin{equation}\label{eq:cauchy}
	\begin{aligned}
		\Delta u &= 0& &\text{ in }\Omega,\\
		u &= 0 &&\text{ on }y=0 \text{ and } y=H,\\
		\frac{\partial u}{\partial x} &= 0 &&\text{ on }x=0,\\
		u &=u_L:= \sin{k\pi\frac{y}{H}}  &&\text{ on }x=0,\\
	\end{aligned}
\end{equation}
where $k\in\mathbb{N}$ is the wavenumber of the signal. As can be observed there is no boundary condition on the right-hand side $\Gamma_R = \{(T,y),y\in(0,H)\}$ whereas there are both Dirichlet and Neumann conditions on the left-hand side. This is a model problem for the non-destructive control of structures in statics, it admits the following solution:
\begin{equation}
	u (x,y)= \sin{k\pi\frac{y}{H}}\cosh(k\pi \frac{x}{H}).
\end{equation}
We observe that the solution on the right-hand side $(x=T)$ explodes for thick domains ($T\nearrow$) or large wavenumbers ($k\nearrow$) as illustrated in Figure~\ref{fig:stekref}. This kind of problem is often qualified as ``severely ill-posed'' \cite{2007_Belgacem_severe}.

\begin{figure}[ht]
	\centering
	\resizebox{0.5\textwidth}{!}{\input{stek_ref.pgf}}\caption{Solution to be identified by Steklov-Poincaré approach observing only the left-hand side (indiscernible oscillation), for $k=3$ and $H=T=1$.}\label{fig:stekref}
\end{figure}

\subsection{Steklov-Poincaré approach}
This approach~\cite{2005_Belgacem_SP} relies on the reformulation of the Cauchy problem in two well-posed problems with one common unknown boundary value $u_R$ on the right-hand side and just one piece of information on the left-hand side.
\begin{equation}
	\begin{aligned}
		\Delta u_D = \Delta u_N &= 0 \qquad&&\text{ in }\Omega,\\
		u_D = u_N &= 0 &&\text{ on }y=0 \text{ and } y=H,\\
		u_D = u_L \text{ and }		\frac{\partial u_N}{\partial x} &= 0 &&\text{ on }x=0,\\
		u_D = u_N& =  u_R&&\text{ on }\Gamma_R. \\
	\end{aligned}
\end{equation}
The index indicates whether Dirichlet or Neumann boundary conditions are considered on the left-hand side. Using classical variational theory, these problems have one solution in $H^{1}(\Omega)$. Using linearity, we can define the Steklov-Poincaré operators $H^{1/2}_{00}(\Gamma_R)\to H^{-1/2}(\Gamma_R)$:
\begin{equation}
	\frac{\partial u_D}{\partial x}= 	\mathcal{S}_D (u_R) - b_D, \qquad \frac{\partial u_N}{\partial x}= 	\mathcal{S}_N (u_R), \qquad\text{on }x=L. 
\end{equation}
$\mathcal{S}_D$ and $\mathcal{S}_N$ are linear operators and $b_D$ is the contribution of the non-zero Dirichlet condition $u_L$.
Solving Cauchy problem~\eqref{eq:cauchy} is then equivalent to finding $u_R$ such that:
\begin{equation}
	(\mathcal{S}_D-\mathcal{S}_N) u_R = b_D.
\end{equation}
The operator $(\mathcal{S}_D-\mathcal{S}_N)$ is compact, meaning that its eigenvalues accumulate near zero. The problem can be discretized using finite element and solved iteratively. We use boldface for the discrete counterpart of the operators and fields introduced above.
In~\cite{2018_Ferrier_Ritz} conjugate gradient with several preconditioners was tested. It was observed that trying to approximate the inverse of the linear operator was a bad idea whereas an efficient preconditioner was provided by the discrete $\mathcal{S}_D$, which in fact corresponded using a Krylov solver to accelerate the KMF stationary iteration~\cite{1991_Kozlov_russe}. Note that $(u_R,\mathcal{S}_D(u_R))=\|u_D\|_{H^1(\Omega)}^2$ so that the KMF preconditioner is in fact a measure of the energy of the field.

\subsection{Direct solvers}
We consider the case $H=T=1$ and $k=3$, the domain is discretized using $40\times 40$ square 4-node Lagrange elements. A Gaussian white noise is added to the signal with signal-to-noise ratio of 10 dB.  The input signal is given in Figure~\ref{fig:signal}. Figure~\ref{fig:direct} presents the identified fields using either a truncated SVD where only singular values larger than ${\varepsilon_\sigma}\sigma_1$ are kept, or a regularized operator $(\mathbf{S}_D-\mathbf{S}_N)+\lambda \mathbf{R})$, $\mathbf{R}$ being either identity $\mathbf{I}$ or $\mathbf{S}_D$. For the record, the five largest eigenvalues of $(\mathbf{S}_D-\mathbf{S}_N)$ are $\{5.8\,10^{-4}$, $2.1\,10^{-6}$, $5.6\,10^{-9}$, $1.2\,10^{-11}$, $2.3\,10^{-14}\}$ while the spectrum of $\mathbf{S}_D$ is contained in the interval $[7.8\ 10^{-3},1.63]$.

\begin{figure}[ht]
	\begin{subfigure}[t]{.49\textwidth}
		\resizebox{0.99\textwidth}{!}{\input{signal.pgf}}\caption{Input signal $u_L$}\label{fig:signal}
	\end{subfigure}
	\begin{subfigure}[t]{.5\textwidth}
		\resizebox{0.99\textwidth}{!}{\input{direct.pgf}}\caption{Identified field $u_R$ with truncated SVD or Tikhonov regularization}\label{fig:direct}
	\end{subfigure}
\end{figure}

For all methods, the straight lines correspond to an identified value of the driving parameter which performs well, while the dotted lines correspond to too strong regularization, and the ``+'' lines correspond to insufficient regularization.
We observe that the Tikhonov regularization by $\mathbf{S}_D$ performs slightly better in the sense that the identification seems to be less sensitive to the value of the regularization parameter. 

\subsection{Preconditioned iterative solvers}

We are now interested in solving the system with conjugate gradient and trying to find efficiently the optimal parameters for the best identification. We compare conjugate gradient without preconditioner $\bM=\bI$, with Jacobi preconditioner $\bM_J = \operatorname{diag}(\bS_D-\bS_N)$ which is usually a good idea for well-posed problems, and the KMF preconditioner $\bM = \bS_D$. The stopping criteria is the one given in~\eqref{eq:stopping2} with $\epsilon=10^{-9}$.

The first column of Figure~\ref{fig:stek:prec} presents the successive approximations during CG (also compared with the reference and with the solution obtained with a truncated SVD with threshold $10^{-12}$). The first three rows correspond to non-regularized systems. It appears clearly that due to the quasi-singularity of the operator, the residual alone is not capable of providing a meaningful stopping criterion: for the first two rows the solvers does not stop at iteration 2 whereas the solution is quite close to the reference. The $\bS_D$ preconditioner appears not to behave significantly differently from the non-preconditioned case, which could be expected as its spectrum does not spread widely. Jacobi preconditioner (third row) behaves poorly, this is mostly due to the fact that the preconditioner is almost singular near the extremities. A fourth row was added which presents the same preconditioner with a slight Tikhonov regularization $\lambda\bS_D$ with $\lambda=10^{-9}$, in that case iterations manage to reduce the error (note that the solver was stopped after 20 iterations).

The second and third columns present the iteration L-curves using either the Euclidean frame $(\|\br\|_2^2,\|\bx_i-\bx_0\|^2_2)$ or the natural frame $(\|\bx_i -\bx\|_\bA^2,\|x_i-x_0\|^2_\bM)$. As could be feared, the Euclidean frame gives hard to exploit curves where it is impossible to define a corner. The natural frame provides much nicer convex curves, and at least in the first two rows, one can clearly see that the last iteration barely reduces the error while strongly increasing the norm of the solution, which suggests selecting the solution at the second to last iteration (which agrees with the ``eye norm'' on the first column).

\begin{figure}
\begin{tabular}{m{.8cm}|c|c|c}
	prec. & Iterations & Euclidean L-curve & Natural L-curve \\\hline
	\rotatebox[origin=lt]{90}{id $\lambda=0$} & 
	\resizebox{0.28\textwidth}{!}{\input{cg_id_0reg_eps9.pgf}}&
	\resizebox{0.28\textwidth}{!}{\input{cg_id_0reg_eps9_lceuclid.pgf}}&
	\resizebox{0.28\textwidth}{!}{\input{cg_id_0reg_eps9_lcnatural.pgf}}\\\hline
	\rotatebox[origin=lt]{90}{KMF $\lambda=0$} & 
	\resizebox{0.28\textwidth}{!}{\input{cg_sd_0reg_eps9.pgf}}&
	\resizebox{0.28\textwidth}{!}{\input{cg_sd_0reg_eps9_lceuclid.pgf}}&
	\resizebox{0.28\textwidth}{!}{\input{cg_Sd_0reg_eps9_lcnatural.pgf}}\\\hline
	\rotatebox[origin=lt]{90}{Jacobi $\lambda=0$} &
	\resizebox{0.28\textwidth}{!}{\input{cg_ja_0reg_eps9.pgf}}&
	\resizebox{0.28\textwidth}{!}{\input{cg_ja_0reg_eps9_lceuclid.pgf}}&
	\resizebox{0.28\textwidth}{!}{\input{cg_ja_0reg_eps9_lcnatural.pgf}}\\\hline
	\rotatebox[origin=lt]{90}{Jacobi $\lambda=10^{-9}$} &
	\resizebox{0.28\textwidth}{!}{\input{cg_ja_9reg_eps9.pgf}}&
	\resizebox{0.28\textwidth}{!}{\input{cg_ja_9reg_eps9_lceuclid.pgf}}&
	\resizebox{0.28\textwidth}{!}{\input{cg_ja_9reg_eps9_lcnatural.pgf}}\\\hline
\end{tabular}\caption{Data completion problem: quality of solution using CG with various preconditioner -- effect of slight Tikhonov regularization on the worst preconditioner.}\label{fig:stek:prec}
\end{figure}

\subsection{Opportunities when preconditioning with Tikhonov operator}
We now consider the KMF operator $\bS_D$ used as both the regularization and the preconditioner: $(\bA+\lambda\bM)\bx = \bb$. After the $m$ CG iterations needed to reach convergence, we postprocess the Ritz values $\Theta^{(m)}$, and the Ritz vectors $\bV^{(m)}$. 
\medskip

Figure~\ref{fig:lambda} presents the approximation of the $\lambda$-L-curve given by the postprocessed Ritz' elements~\eqref{eq:l-curves2} from one iterative computation with $\epsilon=10^{-9}$ and $\lambda=10^{-9}$. It is compared to the L-curve obtained from direct solves with various $\lambda\in[10^{-6,10^{-12}}]$, using the fact that in this academic case the reference $\bx$ is known. Note that because for the iterative solver the actual error is known up to the additive factor $\|x-x_0\|^2_{\bA}$, the curves are in fact translated so that their vertical asymptote is aligned with the $y$-axis.

We see that the L-curve postprocessed at negligible cost after the iterative solution provides an excellent estimation of the actual L-curve which would not be computable in real cases.

\begin{figure}[ht]\centering\resizebox{0.6\textwidth}{!}{\input{lcurves.pgf}}
	\caption{L-curve obtained with direct solves vs L-curve estimated by Ritz' approximation \eqref{eq:l-curves2} after one iterative resolution with $\lambda=10^{-9}$ materialized by the dot $(\epsilon=10^{-9})$. The vertical asymptotes are aligned with the $y$-axis. }\label{fig:lambda}
\end{figure}

\medskip
Figure~ref{fig:picardstek} presents the Picard plot obtained after the iterative resolution with  $\lambda=10^{-9}$ and $\varepsilon=10^{-12}$ (for this simple case, higher precision was required in order to have more points to plot). First, we can observe the spectrum of $\bA$ estimated by the Ritz values. Thanks to~\eqref{eq:specritz}, the regularization is a horizontal line, and we can directly see how the regularization impacts the spectrum: in this case, the two smallest eigenvalues are rectified.

A second interest of Picard plot is to compare the spectrum with the contribution of the right-hand side $|r_j^{(m)}|$. Clearly, beside the first two components, the right-hand side decreases slower than the Ritz' values. This indicates that the system is unstable as the effect of the contributions to the solution gets increasingly larger. 

These two pieces of information observable in Picard plots enable informed choices of $\lambda$ and of the truncation index $i$ in the Ritz reconstruction $\tilde{\bx}^{(m)}_{\lambda,i}$.

\begin{figure}[ht]\centering\resizebox{0.6\textwidth}{!}{\input{picardstek.pgf}}
	\caption{Picard plot processed after the iterative resolution with $\lambda=10^{-9}$ and $\varepsilon=10^{-12}$.}\label{fig:picardstek}
\end{figure}

\section{Augmented preconditioned conjugate gradient for sequences of regularized systems}\label{sec:APCG}
We now consider of nonlinear problems whose solution can be obtained by solving a sequence of linear systems with constant left-hand side and varying right-hand side:
\begin{equation}
	\text{For given }\lambda, \text{ at outer iteration }k \left\{\begin{aligned}
&\text{Compute right-hand side }\bb_\lambda^k(\bx^k_\lambda) \\
&\text{Solve } 
\bA_\lambda \bdelta_\lambda^k = \bb_\lambda^k\\
&\text{Update } \bx^{k+1}_\lambda = \bx^{k}_\lambda + \bdelta_\lambda^k 
	\end{aligned}\right.
\end{equation}

To fully benefit the repetitive nature of the solves, we introduce augmentation in the conjugate gradient. Augmentation is a technique where the search space is split in two subspaces: the augmentation space, defined as the range of the given $n\times n_C$ full-rank matrix $\bC$, and its $\bA$-orthogonal complementary subspace. The part of the solution in the augmentation space is obtained at the initialization whereas the complement is searched for iteratively while maintaining the residual orthogonal to the augmentation. 

We introduce the augmented  Krylov subspace
$\sK_i(\bM^{-1}\bA,\bC,\bM^{-1}\br_0)$ \cite{SAAD:1997:DAK}:
\begin{equation}
	\sK_i(\bM^{-1}\bA,\bC,\bM^{-1}\br_0) = \\
	\operatorname{span}\left(\bM^{-1}\br_0,\ldots,(\bM^{-1}\bA)^{(i-1)}\bM^{-1}\br_0\right)\oplus\operatorname{Range}\left(\bC\right)
\end{equation}
Given an arbitrary initialization $\bx_{00}$ and associated residual $\br_{00}=\bb-\bA\bx_{00}$, the $i_{th}$ iteration can be defined as:
\begin{equation}\label{eq:Akrysol_principle}
	\left\{\begin{array}{ll}
		\textrm{find } &\bx_i \in \bx_{00}+\sK_i(\bM^{-1}\bA,\bC,\bM^{-1}\br_{00})\\
		\textrm{such that }&\br_i \perp \sK_i(\bM^{-1}\bA,\bC,\bM^{-1}\br_{00})
	\end{array}\right.
\end{equation}
This iteration is achieved by Algorithm~\ref{alg:APCG} where the augmentation is a managed by the correction of the initialization (in order to obtain $\bx_0$) and the projector $\bP$ on $\operatorname{ker}(\bC^T\bA)$, which together ensure that the residual remains orthogonal to $\operatorname{Range}(\bC)$ \cite{DOSTAL:1988:CGM}. 

\begin{algorithm}[ht]
	\caption{Augmented Preconditioned Conjugate Gradient}
	\label{alg:APCG}
	\begin{algorithmic}
		\State \color{blue} $\bx_{00}$ and $\bC$ given 
		\State $\bP = \bI - \bC(\bC^T\bA\bC)^{-1}\bC^T\bA$
		\State $\bx_0=\bP \bx_{00}+\bC(\bC^T\bA\bC)^{-1}\bC^T\bb= \bx_{00}+\bC(\bC^T\bA\bC)^{-1}\bC^T\br_{00}$ \color{black}
		\State $\br_0=  \bb - \bA \bx_{0}\color{blue}=\bP^T\br_{00}$\color{black}
		\State {${\bz_{0}} = \color{blue}\bP\color{black}\bM^{-1}   \br_{0}$, $\bw_0=\bz_0$}
		\State $\gamma_0 = (\bz_{0}^T \br_{0})$
		\For{$i = 0,\,1,\, \dots,m$ (convergence)}
		\State  $\bq_i = \bA{\bw_{i}}$
		\State  $\delta_i = (\bw_i^T \bq_{i})$,   $\alpha_i = \delta_i^{-1}\gamma_i$ 
		\State   $\bx_{i+1} = \bx_{i}+\bw_i \alpha_i $
		\State   $\br_{i+1} = \br_i- \bq_i \alpha_i$
		\State   {${\bz_{i+1}} =\color{blue}\bP\color{black}\bM^{-1} \br_{i+1} $}
		\State $\gamma_{i+1} = (\bz_{i+1}^T \br_{i+1})$
		\State  $\beta_{i}= \gamma_i^{-1}\gamma_{i+1}$
		\State $\bw_{i+1} = \bz_{i+1} + \bw_i\beta_{i} $
		\EndFor
	\end{algorithmic}
\end{algorithm}

\subsection{Roles of augmentation}

One particular use of augmentation is to handle symmetric positive semi-definite preconditioner $\bM$: by ensuring $\operatorname{ker}(\bM)\subset\operatorname{Range}(\bC)$ we make sure that the (pseudo-)inverse of $\bM$ is well-defined for the residuals it is applied to. Let $\textbf{\bC}_0$ be a basis of $\operatorname{ker}(\bM)$.\medskip

An additional use of augmentation is to reuse numerical information generated during one solve to accelerate the following. Indeed, augmentation comes with  optimized block operations that make augmenting by one vector much cheaper than one iteration.

In our case, the first solve (in $m$ iterations) will serve to postprocess the Ritz basis $\bV_{m}$ and the next solves will be augmented by the concatenation $\bC=\left[\bC_0,\bV_{\tilde{m}}\right]$ where $\tilde{m}\leqslant m$ indicates that only a selection of the first Ritz vectors is used. Indeed, storing too many vectors can result in excessive memory usage, moreover, the first Ritz vectors are often better converged and improve convergence further than the last  \cite{2013_Gosselet_reuse}.
\medskip

Moreover, Ritz vectors possess two advantages. Firstly, the product $\bA\bV_m$ which is required during augmentation can be obtained at low computational cost using the formula:
\begin{equation}\label{eq:AVm1}
	\begin{aligned}
	\bA\bV_m&=\bA\hbZ_m\bXi_m, \\
	\text{and }\bA\hbz_{j+1} &= (-1)^{j+1}(\bq_{j+1}-\beta_j\bq_j)/\sqrt{\gamma_{j+1}}.
	\end{aligned}
\end{equation} 
Secondly, using normalization:
\begin{equation}\label{eq:AVm2}
	\bV_m \leftarrow \bV_m\bTheta_m^{-1/2} \text{ leads to }\bV_m^T\bA\bV_m=\bI.
\end{equation} 

\subsection{Postprocessing of multiple solutions for various regularization}
The nonlinearity makes it impossible to compute a full analytical L-curve as with formula~\eqref{eq:l-curves2}. Nevertheless, we can choose a family of $P$ regularization coefficients $(\lambda_p)_{p\leqslant P}$ for which we can obtain a good approximation at low cost.
\begin{equation}\label{eq:multiplelambdanl}
\text{At outer iteration }k \left\{\begin{aligned}
		&\text{Compute right-hand sides }\bb_{\lambda_p}^k(\bx^k_{\lambda_p}), \text{ for all }  0\leqslant p\leqslant P\\
		&\text{Solve } 
		\bA_{\lambda_0} \bdelta^k_{\lambda_0} = \bb^k_{\lambda_0} \text{ extract Ritz basis }\bV_{m_k},\\
		& \text{Estimate }\tilde{\bdelta}_{\lambda_p}^{k,(m_k)}\text{ for all } 0< p\leqslant P\text{ using }\eqref{eq:ritzapprox},\\
		&\text{Update } \bx^{k+1}_{\lambda_p} = \bx^{k}_{\lambda_p} + \tilde{\bdelta}_{\lambda_p}^{k,(m_k)},\text{ for all }  0\leqslant p\leqslant P.
	\end{aligned}\right.
\end{equation}
In words, at each outer iteration, only the system associated with $\lambda_0$ is solved with APCG whereas the solutions for the remaining coefficients $(\lambda_p)_{0<p}$ are obtained thanks to Ritz' approximation.

\section{Assessment in the nonlinear case: recovery of the optical flow}\label{sec:OF}
\subsection{Optical flow in a nutshell}\label{ssec:of}

The optical flow is a digital image correlation technique which aims at estimating the displacement field between two images at the scale of the pixel. Contrarily to very popular approaches in solid mechanics inspired by the Finite Element Method \cite{besnard:2006}, it does not rely on a mesh and on shape functions to approximate the displacement field. Given a sequence of two images $(I_1,I_2)$, viewed as $N\times M$ arrays of gray level pixels (in the discrete segment $\{0,1,\ldots,G_{max}\}$), it directly aims at finding the transformation $\phi=(\phi_x,\phi_y)$ such that $I_1-I_2\circ \phi=0$. 
Note that we use interpolation between pixels so that the images can be defined on the rectangle $[0,N]\times[0,M]\subset\R^2$ with values in the continuous segment $[0,G_{max}]\subset \R$ and the displacement $u:=(\phi-\mathcal{I})$ can take non-integer values ($\mathcal{I}$ is the identity operator). It is even common to obtain precision below one tenth of a pixel. 
In order to gain flexibility, and adapt to unavoidable noisy measurements which make the zero unachievable, the problem is better rephrased in terms of the minimization of the ``image energy'' $E_I$:
\begin{equation}
	E_I^2 = \frac{1}{2}\left\|I_1-I_2\circ \phi\right\|^2,\quad\text{where } \|I\|^2=\sum_{\substack{0\leqslant i<N\\0\leqslant j<M}} I(i,j)^2.
\end{equation}
Even under that form the problem is not well-posed, would it only be because there are two times more unknowns than equations. A solution to recover a well-posed problem is to enforce regularity to the displacement field. A penalty term related to the gradient is then introduced:
\begin{equation}
	E^2 = E_I^2 + \frac{\lambda}{2} \|\nabla u\|^2,
\end{equation}
where we kept the Euclidean norm notation for $\|\nabla u\|^2:=\|\partial_x u_x\|^2 + \|\partial_y u_y\|^2 + \|\partial_x u_y\|^2+\|\partial_y u_x\|^2$. $\lambda$ is a weight that needs to be tuned in order to balance the contributions of the image energy and of the regularization. Note that the quadratic minimization framework employed here is not necessarily the most relevant in terms of quality of the identified fields~\cite{Sun2010,chabib:hal-04251608} and that iterative methods were also developed for more sophisticated metrics~\cite{Osher2005}.  \medskip

A modified Gauss-Newton approach is used to minimize the energy \cite{passieux19gaussnewton} which we combine with a pyramidal approach in order to provide meaningful initializations (a sequence of reduced systems is defined, and starting from the coarsest the solution obtained at one level is extrapolated on the next level to define a sound initialization). 
For a given level of the pyramid, starting from a guess $u$, the update $u+du$ is computed by solving the system:
\begin{equation}\label{eq:of}
	( \bA + \lambda \bM )\bx = \bb_{\bA} + \lambda \bb_{\bM},
\end{equation}
with
\begin{equation}\label{eq:of2}
	\begin{aligned}
		\bA &= \begin{pmatrix}
			\bJ_{x} &\\ & \bJ_{y} 
		\end{pmatrix} \begin{pmatrix}
			\bI& \bI\\\bI&\bI 
		\end{pmatrix}\begin{pmatrix}
			\bJ_{x} &\\ & \bJ_{y} 
		\end{pmatrix},\qquad
		\bM = \begin{pmatrix}
			\bDel & \\ & \bDel
		\end{pmatrix}\\ \bx &= \begin{pmatrix}
			\operatorname{vec}(d u_x)\\  \operatorname{vec}(d u_y)
		\end{pmatrix},\qquad
		\bb_{\bA}  = \begin{pmatrix}\operatorname{vec}((I_1-I_2\circ \phi)J_x)\\\operatorname{vec}((I_1-I_2\circ \phi)J_y)	\end{pmatrix},\qquad \bb_{\bM} = \begin{pmatrix}
			\operatorname{vec}(\Delta u_x) \\ 		\operatorname{vec}(\Delta u_y)
		\end{pmatrix}.
	\end{aligned}	
\end{equation}
The $\operatorname{vec}$ operator converts images to vectors ($N\times M$ array to $NM$ vector). For $z\in\{x,y\}$, $J_z$ is the $z$ component of the gradient of $I_1$, $\Delta u_z$ is the (scalar) Laplace operator applied to $u_z$. $\bJ_z$ is the $NM$ diagonal operator containing the values of the gradient $J_z$. $\bI$ and $\bDel$ are respectively the $NM$ identity matrix and the $NM$ matrix version of Laplace operator (with Neumann boundary conditions). All the operators are in fact obtained by discrete difference on the image. Note that the gradient of $I_1$ is used to approximate the current Jacobian. As commonly done in image treatment, a median filter is applied to all the computed increments in order to remove outliers caused by the imperfect speckle.

\begin{figure}[ht]\centering
	\includegraphics[width=.3\textwidth]{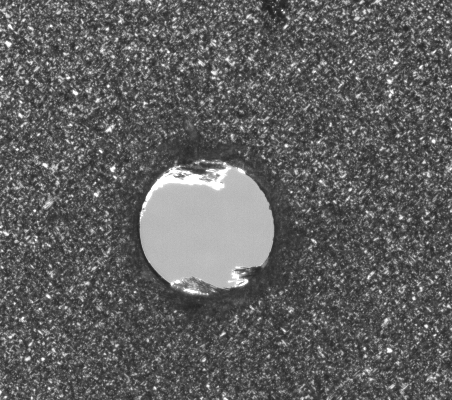}\caption{Speckle of the test specimen.}\label{fig:speckle}
\end{figure}

The proposed test case is a holed composite plate in traction, with a 45$^\circ$ crack to be identified at the bottom of the hole. The speckle in the initial configuration is shown in Figure~\ref{fig:speckle}. For the illustrations, we  present the strain field obtained by deriving the computed displacement. Strain is indeed the mechanical quantity of interest to identify the crack. Only the $xx$ component is given as other components do not provide any further information. 
To quantify the bad conditioning, the non-zero eigenvalues of $\bA$ are in the interval $[10^{-6},10^2]$.

The test case as well as the method (with a simpler linear solver) are presented in~\cite{2023_chabib:gcpu}. For simplicity, we focus on the last nonlinear system to be solved, associated with the full image. Anyhow, the initialization of this system, resulting from the pyramidal approach, was impacted by the choice of the regularization.

\subsection{Matrix free implementation}\label{sec:NLPCG}
We consider the solution to system~(\ref{eq:of}, \ref{eq:of2}) with preconditioned conjugate gradient.

Due to the rectangular shape of the images, there exists an extremely cheap way to solve the regularization matrix, which is a Laplace operator, using Fast Fourier transform or more precisely discrete cosine transform, see Appendix~\ref{sec:invlaplace}. 

It is extremely simple to work with $\bA$ and $\bM$ without assembling them, one only needs to compute and store the two $N\times M$ images  $(\bJ_x,\bJ_y)$ and use Hadamard product and Laplace function when computing matrix-vector multiplication.
\medskip

The system is of dimension $2MN$. As said earlier, $\bA$ is strongly deficient since its rank is at most $MN$, a first part of its kernel has the following basis:
\begin{equation}
	\operatorname{span}\begin{pmatrix}
		\bJ_y\\-\bJ_x\end{pmatrix} \subset \operatorname{ker}(\bA).
\end{equation} 
The rest of the spectrum is easy to compute since:
\begin{equation}
	\bA \begin{pmatrix}
		\bJ_x\\\bJ_y
	\end{pmatrix} = \begin{pmatrix}
		\bJ_x\\\bJ_y
	\end{pmatrix}\left(\bJ_x^2+\bJ_y^2\right).
\end{equation}
The other eigenvalues thus correspond to the square of the norm of the gradient of the image. Pixels where the gradient is zero (bad speckles) are also associated with zero eigenvalues.

$\bM$ is also rank deficient, the dimension of its kernel is 2, a basis of its null space is well known:
\begin{equation}
	\operatorname{ker}(\bM)= \operatorname{span}\begin{pmatrix}
		\bone&0\\0&\bone
	\end{pmatrix},
\end{equation} 
where $\bone$ is the vector filled with 1: the kernel of the scalar Laplace operator consists of constant functions. In fact a more efficient basis can be computed at a very low cost:
\begin{equation}\label{eq:kerbasis}
	\bC_0 = \begin{pmatrix}
	\displaystyle	\frac{ 1}{\sqrt{s_{xx}}}\bone& \displaystyle-\frac{s_{xy}s_b}{
		{s_{xx}}}\bone \\ 0 & \displaystyle s_b\bone
	\end{pmatrix} \text{with}\left\{\scriptstyle\begin{aligned}
		&	s_{xx} = \bone^T\bJ_x^2\bone =\sum_i j_{x,ii}^2\\
		&		s_{xy} = \bone^T\bJ_x\bJ_y\bone =\sum_i j_{x,ii}j_{y,ii}\\
		&		s_{yy} = \bone^T\bJ_y^2\bone=\sum_i j_{y,ii}^2 \\
		&   s_b = 1/\sqrt{s_{yy}-s_{xy}^2/s_{xx}}
	\end{aligned}\right.
\end{equation}
This matrix is used as an augmentation in order to make sure that we only work in a space orthogonal to the kernel of $\bM$. 
It has the advantage to make the matrix  $(\bC_0^T(\bA+\lambda\bM)\bC_0)=(\bC_0^T\bA\bC_0)=\bI$ for any $\lambda$.

\subsection{Quality of the preconditioner}
We first wish to verify that preconditioning by regularization actually leads to better enforcement of the regularity. In Table~\ref{tab:compa_prec}, we can qualitatively compare the classical Jacobi approach of preconditioning by the diagonal of the operator $\bM_{jac}^{-1}=\operatorname{diag}(\bA_\lambda)^{-1}$ and the proposed preconditioning by regularization. The increased regularity is particularly visible for low weight $\lambda$ and low precision $\varepsilon$ of the linear solver.

Preconditioning by the regularization operator thus makes it possible to make meaningful computations with low weight in the regularization and to solve with less precision, hence with fewer iterations. Nevertheless, one has to mention that our preconditioner is computationally more expensive per iteration than the diagonal one.\medskip

\begin{table}[ht]\centering
	\begin{tabular}{|C{1cm}|C{1cm}|C{5.cm}|C{5.cm}|}\hline
		$\lambda$ & $\varepsilon$ & Diagonal Prec. & Regularization Prec.\\\hline
		low & low & \includegraphics[width=\linewidth]{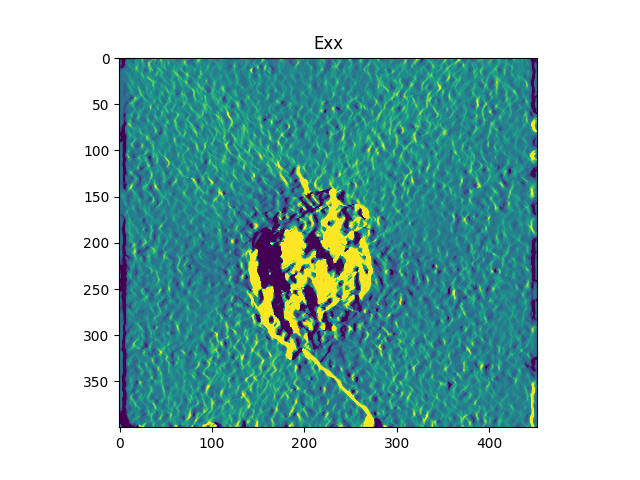} & \includegraphics[width=\linewidth]{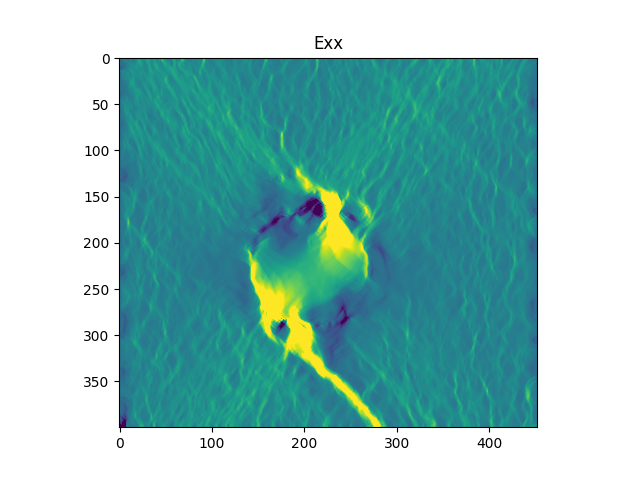}\\\hline
		low & high & \includegraphics[width=\linewidth]{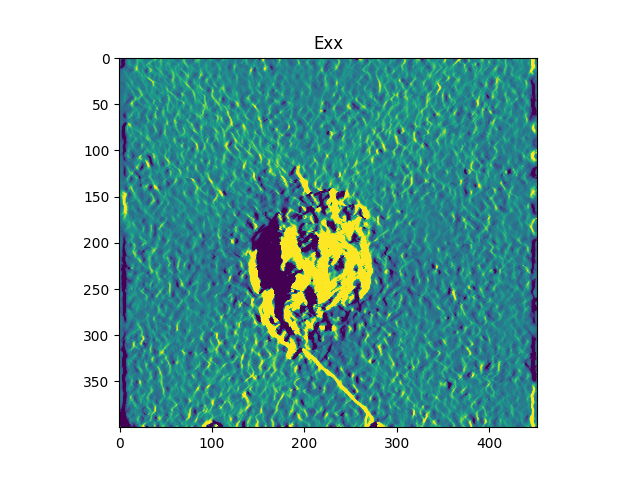} & \includegraphics[width=\linewidth]{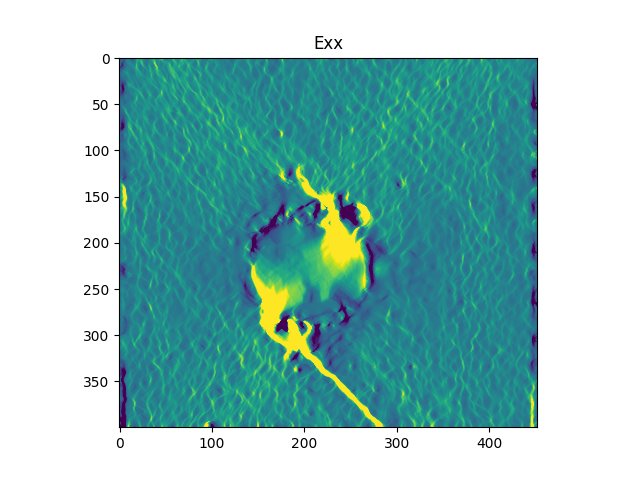} \\\hline
		high & low  & \includegraphics[width=\linewidth]{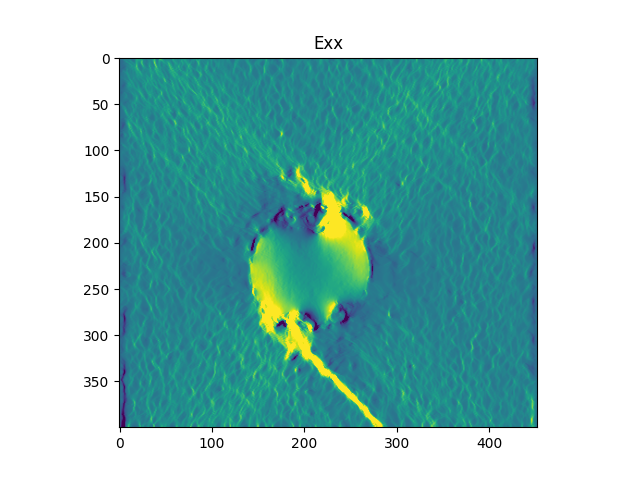} & \includegraphics[width=\linewidth]{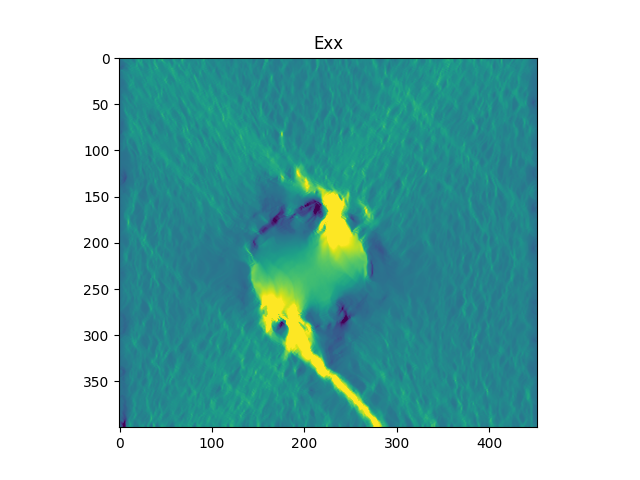} \\\hline
		high & high  & \includegraphics[width=\linewidth]{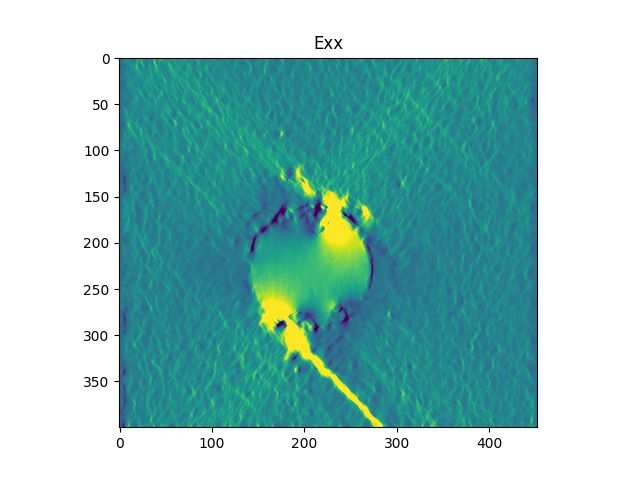} & \includegraphics[width=\linewidth]{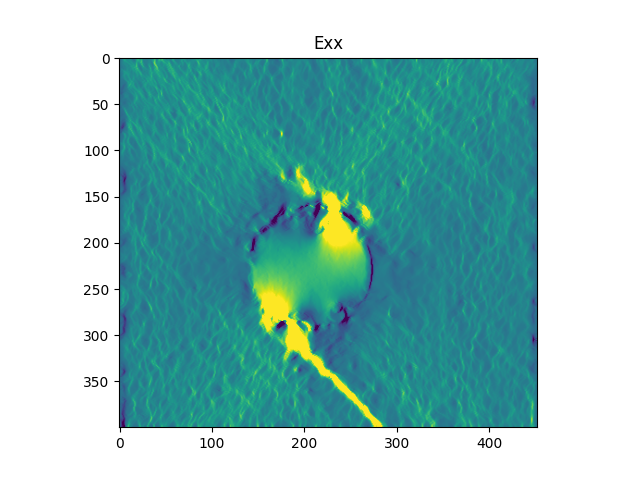} \\\hline
	\end{tabular}\caption{$\epsilon_{xx}$ strain field (range = mean value $\pm$ 3 st.dev.). Comparison of the effect of preconditioning by diagonal (simple approach) vs by regularization, for different weights $\lambda\in\{1,1000\}$ and linear solver precision $\varepsilon\in\{10^{-2},10^{-3}\}$.}\label{tab:compa_prec}
\end{table}

An interesting scenario unfolds. The regularization preconditioner promotes low frequency corrections. Indeed, it is associated with a fully populated matrix (never actually computed) and the search directions have naturally large wavelength. As iterations progress, higher frequency modes emerge introducing more and more details and irregularity. On the contrary, the Jacobi preconditioner is diagonal, and it naturally encourages (independent) details, only iterations make it possible to reveal the structure between neighboring pixels.

\subsection{Ritz filtering}

\begin{figure}[ht]
	\null\hfill
	\begin{subfigure}{.48\textwidth}
		\centering\includegraphics[width=.95\textwidth]{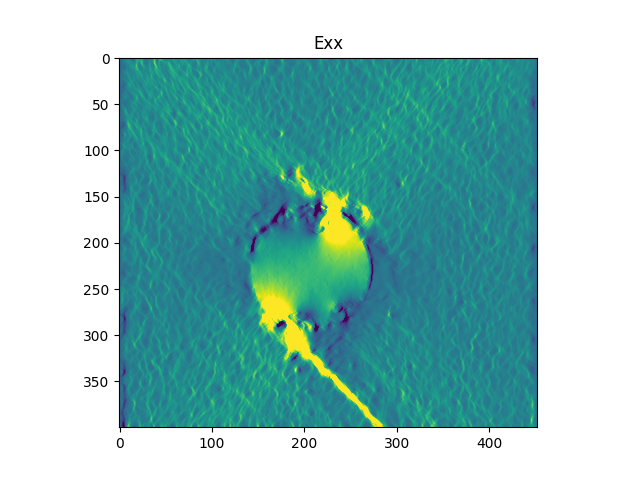}\caption{$\lambda=1000$}
	\end{subfigure}\hfill
	\begin{subfigure}{.48\textwidth}
		\centering\includegraphics[width=.95\textwidth]{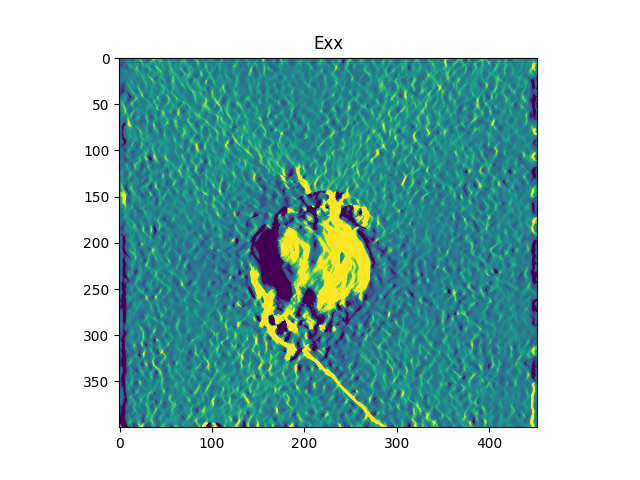}\caption{$\lambda=10$}
	\end{subfigure}\hfill\null\caption{Identified $\epsilon_{xx}$ field, for different regularization intensity $\lambda$, with $\varepsilon=10^{-5}$.}\label{fig:Exx}
\end{figure}

We analyze the solving process for high ($\lambda=1000$) and low ($\lambda=10$) levels of regularization. We use the second stopping criterion of Equation~\eqref{eq:stopping} with $\varepsilon=10^{-5}$, which corresponds to a rather high level of convergence. The identified strain field are given in Figure~\ref{fig:Exx}.

\begin{figure}[ht]
	\null\hfill
	\begin{subfigure}{.48\textwidth}
		\centering\includegraphics[width=.95\textwidth]{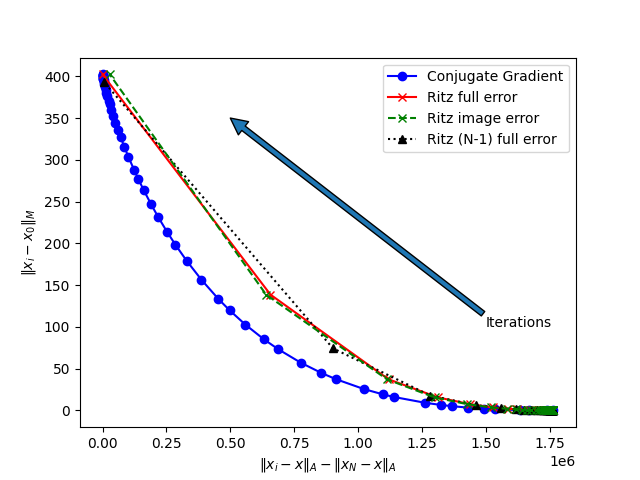}\caption{$\lambda=1000$, 56 iterations}
	\end{subfigure}\hfill
	\begin{subfigure}{.48\textwidth}
		\centering\includegraphics[width=.95\textwidth]{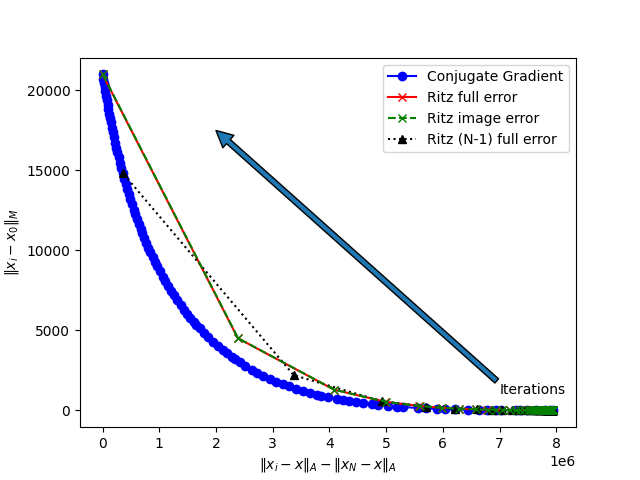}\caption{$\lambda=10$, 104 iterations}
	\end{subfigure}\hfill\null\caption{L-curves (same $\varepsilon=10^{-5}$) compared with Ritz post-treatment, for different regularization intensity $\lambda$. Note that different scales are used on the plots. }\label{fig:Lcurve}
\end{figure}

We analyze the convergence in terms of compromise between the decrease of the error and the increase of the norm of the gradient of the solution which stems from the oscillations in the identified fields. Figure~\ref{fig:Lcurve} presents two L-curves associated with high and low regularization. We use the natural CG-norms, please note that the position of the 0-abscissa is conventional because $\|\bx_0-\bx\|_{\bA}$ is unknown. The L-curves of the CG iterations (dotted lines) have similar shapes, like pieces of hyperbola. Due to the difference of magnitude, different scales had to be used: the error decreases four times less when the high regularization is used, and the norm of the solution remains 50 times smaller.

In order to better understand the convergence, we conduct a Ritz analysis.  For $\lambda=1000$, the convergence is attained in $m=56$ iterations and as many Ritz vectors are computed. Table~\ref{tab:compa_ritz} presents a selection of these modes, sorted in decreasing order of Ritz value. The Ritz vectors resemble vibration modes with increasing number of anti-nodes. The first vectors are so regular that the hole is barely visible. The crack is only visible on the latest modes. This is bad (but logical) news because these are the most difficult modes to converge, thus they are probably bad approximations of actual eigenmodes, and the crucial mechanical information they carry is difficult to reuse.

\begin{table}[ht]\centering
	\begin{tabular}{|C{3.5cm}|C{3.5cm}|C{3.5cm}|C{3.5cm}|}\hline
		$\theta_0=11\,276\,321$& $\theta_{17}=966\,277$ &  $\theta_{27}=669\,871$& $\theta_{39}=256\,321$ \\
		\includegraphics[width=\linewidth]{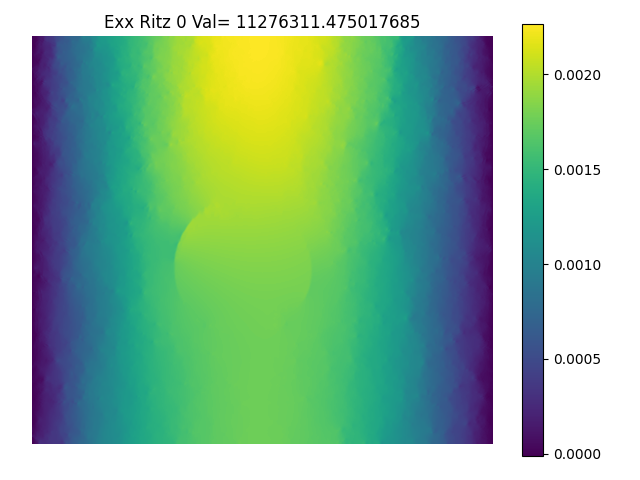} & \includegraphics[width=\linewidth]{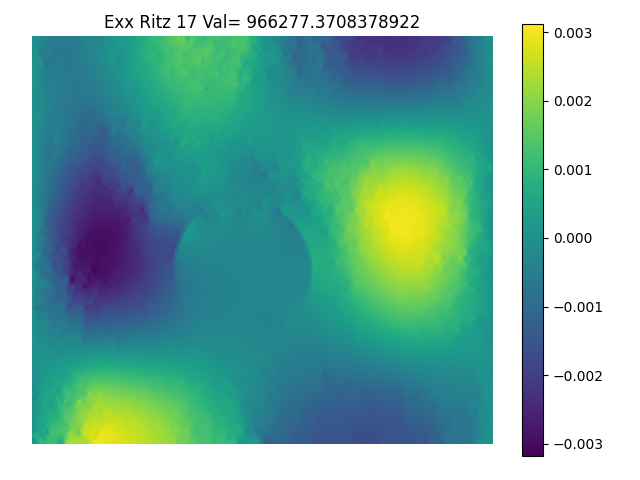}&
		\includegraphics[width=\linewidth]{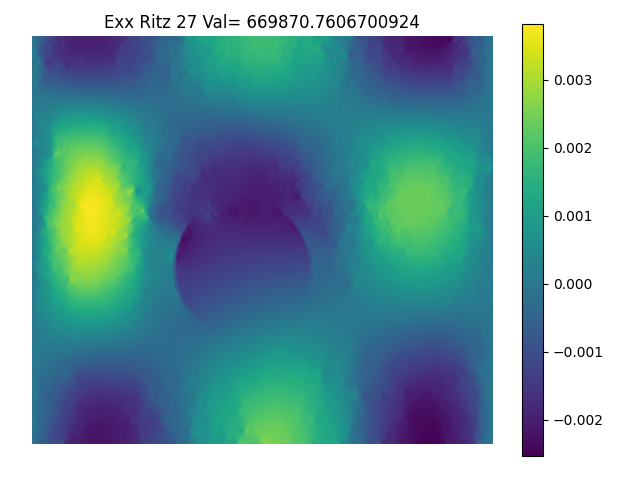} & \includegraphics[width=\linewidth]{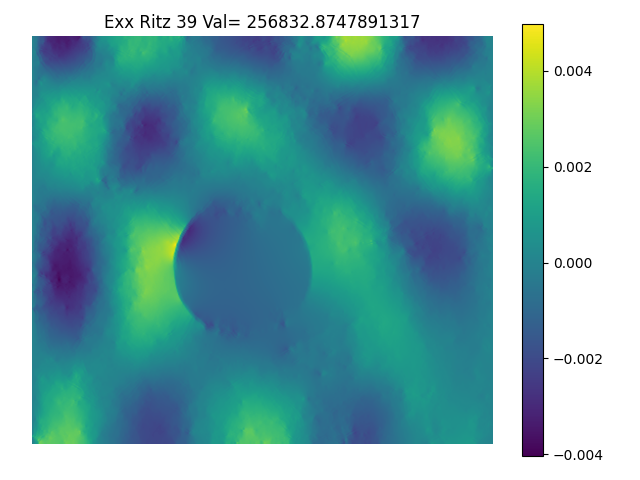} \\\hline
		$\theta_{47}=77\,671$& $\theta_{51}=21\,520$ &  $\theta_{54}=2\,560$& $\theta_{55}=493$ \\
		\includegraphics[width=\linewidth]{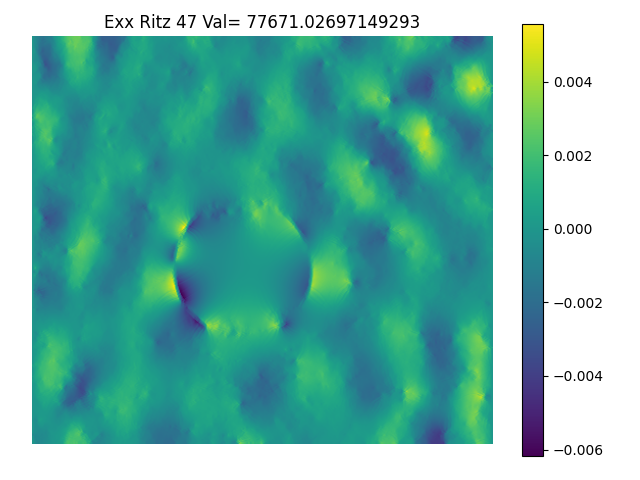} & \includegraphics[width=\linewidth]{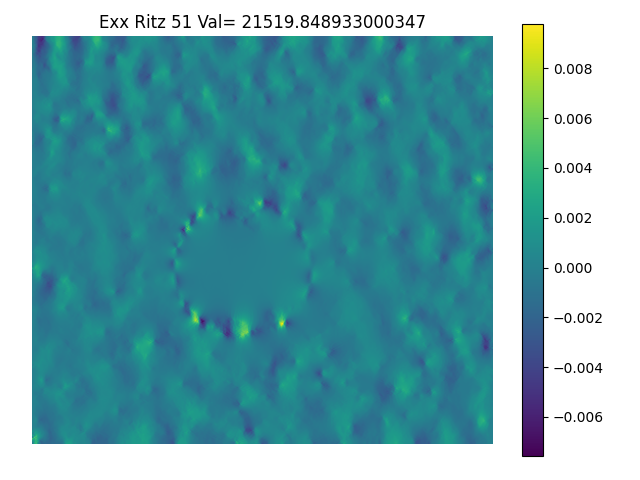} &
		\includegraphics[width=\linewidth]{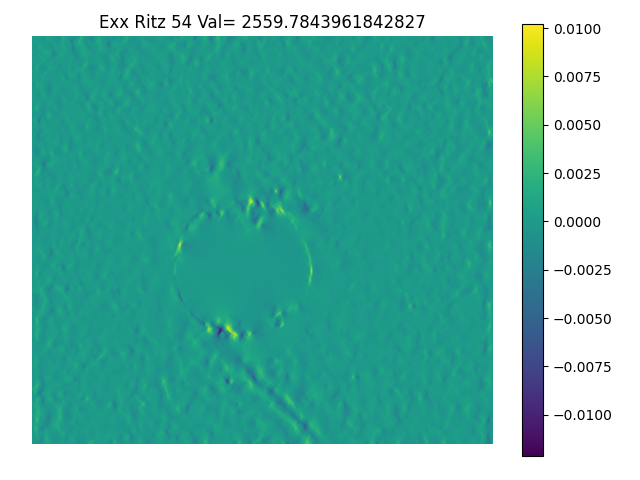} & \includegraphics[width=\linewidth]{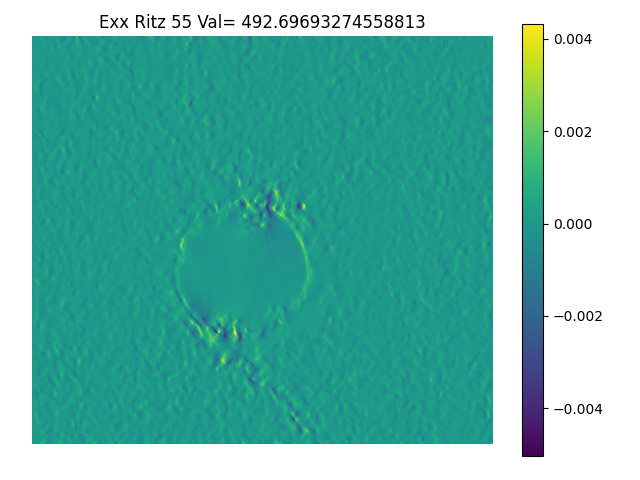} \\\hline
	\end{tabular}\caption{$\epsilon_{xx}$ strain field for 8 Ritz vectors out of the 56 computed ($\varepsilon=10^{-5}$, $\lambda=1000$).}\label{tab:compa_ritz}
\end{table}

We use formula \eqref{eq:ritzapprox} for the \textit{a posteriori} filtering of the solution based on Ritz vectors.
We present the L-curves in terms of modes included in the reconstruction. We show the curves in terms of full error and only taking into account the image error $(\bv_i^T\br_{\bA,0})$ --- they are almost  overlaid on each other, a slight discrepancy only appears for high regularization. The shape of the Ritz L-curves corresponds to most of the modes (the highest) only slightly decreasing the error and almost not changing the norm, only the last modes, which contain the crack information, are associated with significant decrease of the error (but of course at the cost of much increased solution norm). More or less, if a corner was to be selected it would correspond to just suppressing the contribution of the last mode.

In order to better understand this behavior, we first analyze the convergence of the Ritz values by comparing the spectrum obtained at the last iteration with the one obtained just one iteration before, like was done in~\cite{2013_Gosselet_reuse} in the case of a well posed problem. It appears that the largest Ritz values were quite well approximated and only the lowest part of the spectrum evolves (remember the Ritz values correspond to the inverse of the slope of the segments in the L-curve). In other words, even though the last iterations seem not to modify the solution much (accumulation of the dots in the upper left part on the CG L-curves), they play an important role in terms of estimation of the lower part of the spectrum, without adding lots of small eigenvalues.

\begin{figure}[ht]
	\null\hfill
	\begin{subfigure}{.48\textwidth}
		\centering\includegraphics[width=.95\textwidth]{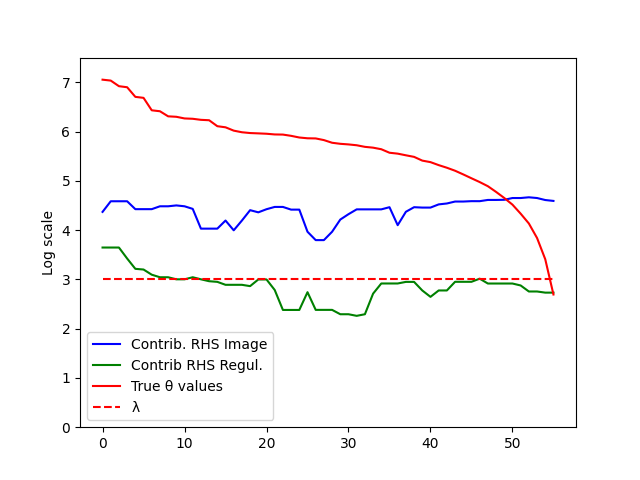}\caption{$\lambda=1000\simeq\min(\theta_j)$}
	\end{subfigure}\hfill
	\begin{subfigure}{.48\textwidth}
		\centering\includegraphics[width=.95\textwidth]{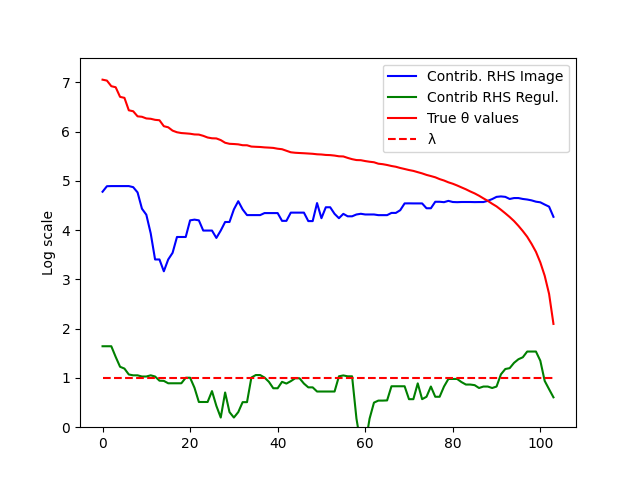}\caption{$\lambda=10\ll(\theta_j)$}
	\end{subfigure}\hfill\null\caption{Spectral analysis of the system for $\varepsilon=10^{-5}$ and different regularization intensity $\lambda$. A 5-width median filter was used to smooth out the contribution curves.}\label{fig:picard}
\end{figure}

To support this analysis, we conduct a Picard's study on Figure~\ref{fig:picard} which shows the distribution of the Ritz values $(\theta^{(m)}_i)$ as well as the decomposition of the right-hand side on the eigenspace $(\bv^{(m)^T}_i\br_{\bA,0})$ and $\lambda(\bv_i^{(m)^T}\br_{\bM,0})$. It is worth recalling that low and highly regularized systems have the same spectrum, except that it is more sampled for the low regularization which requires two times more iterations to converge. The Ritz values are slowly decreasing and only the last 10\% really decay, the low regularization is not associated with an overpopulation of the lowest part of the spectrum. What stands out is the fact that the right-hand side contributes almost equally on all modes (at least it does not decrease for larger Ritz values). Picard's theory thus suggests that we should stop the reconstruction when the Ritz values start to decay. This is not possible in our case since the crack is mostly represented in this part of the spectrum.

By the way, Figure~\ref{fig:picard} permits to compare the smallest Ritz value $\theta_{\min}$ with the regularization parameter $\lambda$. The case that we called ``low regularization'' corresponds to $\lambda$ being negligible with respect to the small Ritz value $\theta^{(m)}_{m}$, and thus only marginally modifying the active Ritz spectrum. On the contrary, the high regularization corresponds to a $\lambda>\theta^{(m)}_{m}$ which means that the lower part of the spectrum of $\bA_\lambda$ is flattened relative to that of $\bA$.

\subsection{Subspace recycling}
Even though it appears that Ritz filtering is difficult to apply to the studied system, we can still benefit from Ritz vectors to accelerate the solution. As a sequence of linear systems with identical matrix has to be solved, it is natural to augment the system with the previously generated Ritz vectors by concatenating  $\bC \leftarrow \begin{pmatrix} 	\bC & \bV_m \end{pmatrix}$. Indeed, augmentation comes with  optimized block operations that make augmenting by one vector much cheaper than one iteration in particular using \eqref{eq:AVm1} and~\eqref{eq:AVm2}.

\begin{table}[ht]\centering
	\begin{tabular}{|l|c|c|c|c|c|c|c|c|c|}\hline
		Aug. &  0   & 10   & 20  & 30  & 40  & 50  & 60  & 70  & max (77) \\\hline
		Iter.&  77  & 64   & 57  & 49  & 44  & 41  & 40  & 40  & 38 \\\hline
		Time (s)   & 11.7 & 10.1 & 8.3 & 7.1 & 6.6 & 6.3 & 6.2 & 6.3 & 6.6 \\\hline
	\end{tabular}\caption{Performance of recycling}\label{tab:perfrec}
\end{table}

Table~\ref{tab:perfrec} illustrates the performance of recycling for the nonlinear system to be solved on the full image at the end of identification. The first linear system, only augmented by the kernel of the preconditioner is solved in 77 iterations. Then a certain portion of the Ritz vectors is used to augment the next 8 linear systems (same matrix, different right-hand sides). As the augmentation results in an excellent initialization, we use a criterion in terms of absolute value of $\|\br_j\|_{\bM}$ to halt the iterations because other comparison as given in Equation~\eqref{eq:stopping} might use an unfair reference. We measure the performance in terms of gain in iterations, and in computational time (measures are conducted on a upper mid-range laptop with Nvidia RTXA2000 graphic card). The gain in terms of iterations is moderate, with best obtained for small augmentation space (at most 1.3 iterations per augmentation vector, for 10 vectors). In terms of time, the optimal is obtained for augmentation space of 80\%-90\% of available vectors, with a global CPU time divided by almost 2 (this time includes all the extra cost associated with computing and using Ritz vectors). This size of subspace agrees with what we observed on the stability of the largest Ritz values in the L-curves plots.

\subsection{Tuning of $\lambda$}
It is often hard to automatize the selection of the regularization intensity $\lambda$. Picard's plots like in Figure~\ref{fig:picard} permit to put $\lambda$ in relation with the spectrum of the preconditioned operator and thus to understand the effect of the regularization in terms of flattened spectrum. Still, the final judge is often the expert's impression of a strain map, and it is convenient to compute maps associated with several $(\lambda_p)$ at low cost. 

As presented in \eqref{eq:multiplelambdanl}, after solving one linearized system regularized by $\lambda_0$, we post-process the solution for any $\lambda_p$ at the simple cost of computing the associated right-hand side $\bb_{\lambda_p}^k(\bx^k_{\lambda_p})$ (which depends on the history of the nonlinear solution for $\lambda_p$), and basic linear algebra operations.

\begin{figure}
	\null\hfill
	\begin{subfigure}{.24\textwidth}\centering
		\hspace{.95\textwidth}
	\end{subfigure}\hfill
	\begin{subfigure}{.24\textwidth}\centering
		\includegraphics[
		width=.95\textwidth]{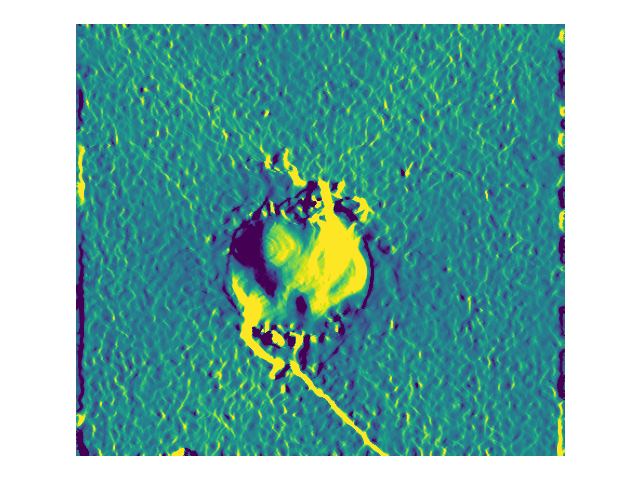}\caption{$\lambda_0=10$}
	\end{subfigure}
	\hfill
	\begin{subfigure}{.24\textwidth}\centering
		\includegraphics[
		width=.95\textwidth]{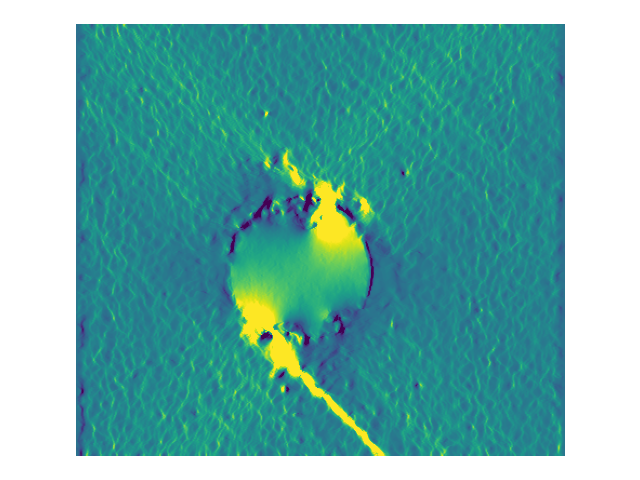}\caption{$\lambda_0=1000$}
	\end{subfigure}
	\hfill
	\begin{subfigure}{.24\textwidth}\centering
		\includegraphics[
		width=.95\textwidth]{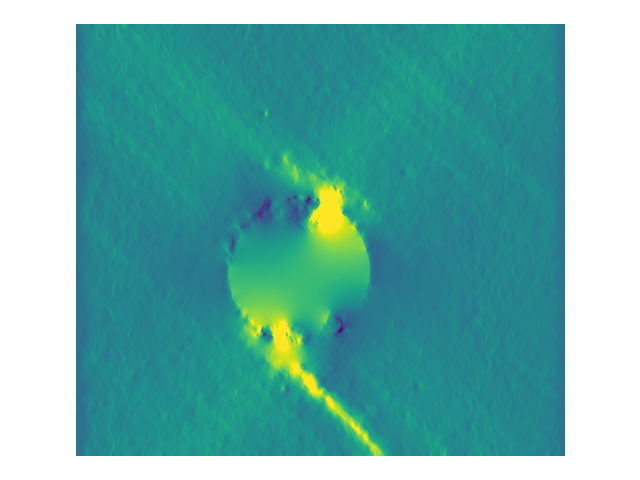}\caption{$\lambda_0=10000$}
	\end{subfigure}
	\null\hfill
	
	\null\hfill
	\begin{subfigure}{.24\textwidth}\centering
		\includegraphics[
		width=.95\textwidth]{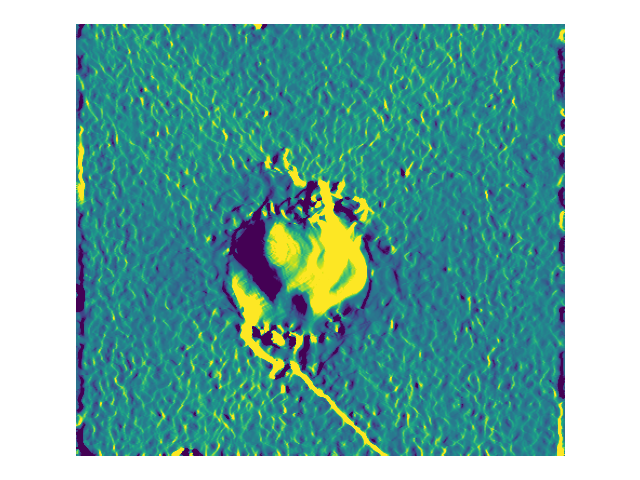}\caption{$\lambda=1$.}
	\end{subfigure}
	\hfill
	\begin{subfigure}{.24\textwidth}\centering
		\includegraphics[
		width=.95\textwidth]{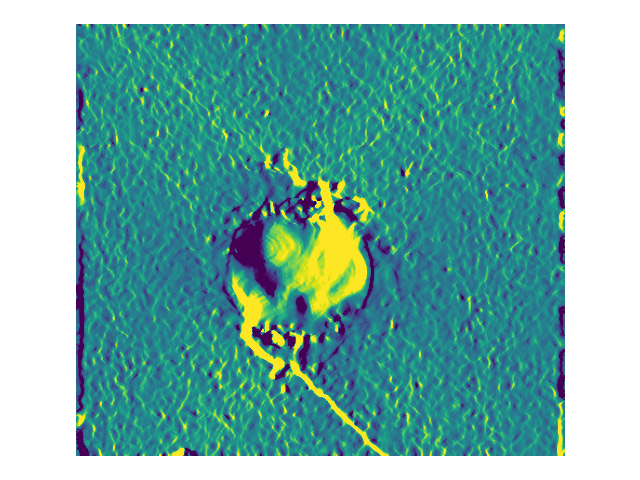}\caption{$\lambda_1=1$  $(\lambda_0=10)$.}
	\end{subfigure}
	\hfill
	\begin{subfigure}{.24\textwidth}\centering
		\includegraphics[
		width=.95\textwidth]{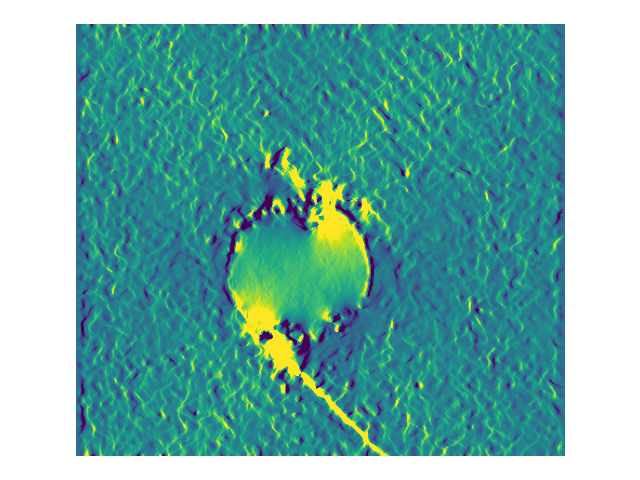}\caption{$\lambda_1=1$  $(\lambda_0=1000)$.}
	\end{subfigure}
	\hfill
	\begin{subfigure}{.24\textwidth}\centering
		\includegraphics[
		width=.95\textwidth]{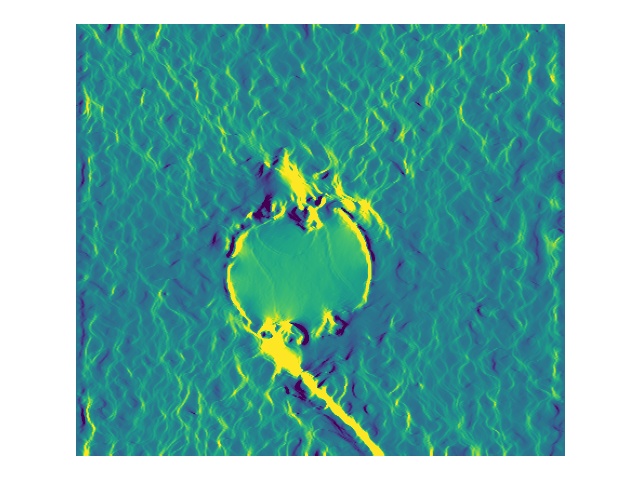}\caption{$\lambda_1=1$  $(\lambda_0=10000)$.}
	\end{subfigure}
	\null\hfill
	\caption{Costless postprocessing for different $(\lambda_i)$: top, initial computation with $\lambda_0\in\{10,100,1000\}$; bottom, direct solution with $\lambda=1$ and solutions deduced with $\lambda_1=1$. $\epsilon_{xx}$ strain field.}\label{fig:tuninglambda}
\end{figure}

Figure~\ref{fig:tuninglambda} presents the solution deduced for $\lambda_1=1$ from initial computations with different  $\lambda_0$ (in $\{10,1000,10000\}$) and $\varepsilon=10^{-4}$. Again, a median filter was applied after the Ritz reconstruction. 
It seems that the Ritz vectors make it possible to postprocess a reasonable solution with up to $\lambda_1\simeq\lambda_0/1000$. The reconstructed strain field appears to be much less smooth than the original computation (with $\lambda_0$) while less noisy than the direct low-regularization computation with $\lambda=1$. If the deduced solution is not fully satisfying, it can still be used as an excellent initialization for a regular computation.


\section{Conclusion}\label{sec:conc}
In this paper, we have studied how preconditioning and Tikhonov regularization could be efficiently combined in an augmented preconditioned conjugate gradient. We have shown that this association makes sense from a physical point of view and it made it possible to combine filtering, recycling of subspaces, and postprocessing of all regularized solutions at negligible cost. This gives a favorable framework to apply criteria like the L-curve or Picard's analysis.

First, the solver was applied to the boundary completion problem. It was verified that regularization was a good preconditioner and that Ritz analysis made it possible to approximate pretty efficiently the L-curve and to build extremely informative Picard plots.

The solver was then applied to a problem of the optical flow reconstruction which introduced the extra difficulty of nonlinearity and the fact that the most important information was buried in the lower part of the spectrum. Preconditioning by the regularization made it possible to start from a regular estimation and enhance details through iterations. Satisfying results were obtained on actual measurements from digital image correlation of a mechanical test. Postprocessed solutions at negligible cost were still relevant for regularization weight $\lambda$ divided by up to $1000$.

An obvious next step for this work is to consider inexact preconditioners, that is to say when the $\bM^{-1}$ matrix in the preconditioning step of the algorithm is only an approximation of the inverse of the regularization matrix $\bM$ in the operator. This would make the method applicable on a much broader class of problems.

\bibliographystyle{plain}
\bibliography{cauchy}

\appendix

\section{Inverse of Laplacian on a rectangle with Neumann boundary condition}\label{sec:invlaplace}

It is well known that plane waves $x\mapsto e^{i\omega\cdot x}$, with $\omega \in\mathbb{R}^2$, form a set of eigenfunctions for the Laplace operator in $\mathbb{R}^2$ with eigenvalues $-\|\omega\|^2$ (using the Euclidean norm). This can be equivalently formulated by saying that the Fourier transform diagonalizes the Laplacian. Hence, the powerful  solution technique (in that case $\omega$ is the variable in the Fourier domain):
\begin{equation}
	\begin{aligned}
		\Delta u +f &= 0\qquad \text{in }\mathbb{R}^2 \\
		u &= \mathcal{F}^{-1} \left( \frac{\mathcal{F}(f)}{\|\omega\|^2}\right)
	\end{aligned}
\end{equation}

What is remarkable is that the eigenvectors are preserved by discretization. For instance, if we consider the classical 5-point stencil on a unit grid:
\begin{equation}
	(\Delta_h u )(x_1,x_2) = u(x_1+1,x_2)+u(x_1-1,x_2)+u(x_1,x_2+1)+u(x_1,x_2-1)-4u(x_1,x_2)
\end{equation}
and one can check that
\begin{equation}
	(\Delta_h e^{i\omega\cdot x})(x_1,x_2) = \underbrace{e^{i (x_1\omega_1+x_2\omega_2)}}_{e^{i\omega\cdot x}}(e^{i\omega_1}+e^{-i\omega_1}+e^{i\omega_2}+e^{-i\omega_2}-4).
\end{equation}

Now, considering a rectangular domain, the boundedness of the domain and the boundary conditions lead to selecting only certain eigenvalues, and eigenvectors are made out of a good combination of plane waves. Consider the unit square $[0,1]^2$, the eigenvalues $\lambda_{n,m}$ and eigenvectors $v_{n,m}$ of the Laplacian with (homogeneous) Neumann boundary conditions are given by:
\begin{equation}
	\begin{aligned}
		v_{n,m}(k,l)& = \cos(\frac{ml\pi}{M})\cos(\frac{nk\pi}{N})\\
		\lambda_{n,m} &= 2 \left(1-\cos(\frac{n\pi}{N})\right) + 2 \left(1-\cos(\frac{m\pi}{M})\right)
	\end{aligned}
\end{equation}
As eigenvectors are cosine functions, the specialization of the Fourier transform to this case takes the name of discrete cosine transform (DCT).

One just needs to take some care of the eigenvalue $\lambda_{0,0}=0$, associated with the constant eigenvector. The classical solution is to work on functions with zero mean value and nullify the constant term in the transformed function.

\begin{figure}[ht]
	\centering\includegraphics[]{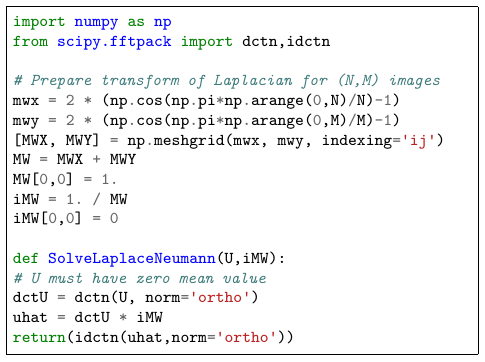}\caption{Inverse discrete Laplacian on a rectangle (Neumann bcs) in python}\label{code:dct}
\end{figure}

Figure~\ref{code:dct} gives the python code for the inverse of the discrete Laplacian on a rectangle with Neumann boundary conditions. This discrete Laplace operator can be directly invoked by the \texttt{laplace()} function  from \texttt{scipy.ndimage} with default arguments (\texttt{border='reflect'}). 

%
%

\end{document}